\documentclass[11pt, twoside,final]{scrartcl}
\usepackage{amssymb, amsfonts, amsmath, mathtools}
\usepackage{theorem}
\usepackage{graphicx}
\graphicspath{{pics/}}
\usepackage{todonotes}
\usepackage[skip=1pt, font=normalsize]{subcaption}
\allowdisplaybreaks

\usepackage{hyperref} 
\hypersetup{plainpages=false, colorlinks, linkcolor=black, citecolor=black, urlcolor=blue,
pdftitle={Nonuniform fast Fourier transforms with nonequispaced spatial and frequency data and fast sinc transforms},
pdfauthor={Melanie Kircheis, Daniel Potts, Manfred Tasche},
pdfstartview={FitBH}}

\theoremstyle{plain}
\newtheorem{theorem}{Theorem}[section]
\newtheorem{corollary}[theorem]{Corollary}
\newtheorem{lemma}[theorem]{Lemma}
\newtheorem{algorithm}[theorem]{Algorithm}
\newtheorem{definition}[theorem]{Definition}
\newtheorem{example}[theorem]{Example}
\newtheorem{remark}[theorem]{Remark}

\newcommand{\bend}{\hspace*{0ex} \hfill \hbox{\vrule height
    1.5ex\vbox{\hrule width 1.4ex \vskip 1.4ex\hrule  width 1.4ex}\vrule
    height 1.5ex}}
\newcommand{\qedsymbol}{\rule{1.5ex}{1.5ex}}

\newenvironment{Lemma}{\goodbreak\begin{lemma}}{\end{lemma}}
\newenvironment{Theorem}{\goodbreak\begin{theorem}}{\end{theorem}}

\newenvironment{Example}{\goodbreak\begin{example}\normalfont}{\bend\end{example}}

\numberwithin{equation}{section}
\numberwithin{table}{section}
\numberwithin{figure}{section}

\renewcommand{\mathbf}[1]{\ensuremath{\boldsymbol{#1}}}
\newcommand{\new}[1]{\textcolor{black}{ #1}}
\newcommand{\I}{\mathcal{I}}

\title{Nonuniform fast Fourier transforms with nonequispaced spatial and frequency data and fast sinc transforms}
\author{
  Melanie Kircheis\footnotemark[4] \and
  Daniel Potts\footnotemark[1]\and
  Manfred Tasche\footnotemark[3]
}

\date{}
\begin{document}
\maketitle

\begin{abstract}
In this paper we study the nonuniform fast Fourier transform with nonequispaced spatial and frequency data (NNFFT) and the fast $\mathrm{sinc}$ transform as its application.
The
computation of NNFFT is mainly based on the nonuniform fast Fourier transform with nonequispaced spatial nodes
and equispaced frequencies (NFFT). The NNFFT employs two compactly supported, continuous window functions.
For fixed nonharmonic bandwidth, we show that the error of the NNFFT with two $\sinh$-type window
functions has an exponential decay with respect to the truncation parameters of the used window functions. As an important application of the NNFFT,
we present the fast $\mathrm{sinc}$ transform. The error of the fast $\mathrm{sinc}$ transform is estimated as well.
\medskip

\emph{Key words}: nonuniform fast Fourier transform, NUFFT, NNFFT, non\-equi\-spaced nodes in space and frequency domain, exponential sums, fast $\mathrm{sinc}$ transform, error estimates, sampling.
\smallskip

AMS \emph{Subject Classifications}:
65T50,
94A12,  
94A20.
\end{abstract}

\footnotetext[4]{melanie.kircheis@math.tu-chemnitz.de, Chemnitz University of
	Technology, Faculty of Mathematics, D--09107 Chemnitz, Germany}

\footnotetext[1]{potts@math.tu-chemnitz.de, Chemnitz University of
	Technology, Faculty of Mathematics, D--09107 Chemnitz, Germany}

\footnotetext[3]{manfred.tasche@uni-rostock.de, University of Rostock, Institute of Mathematics, D--18051 Rostock, Germany}

\section{Introduction}

The \emph{discrete Fourier transform} (DFT) can easily be generalized to arbitrary nodes in the space domain as well as in the frequency domain (see \cite{DuRo93, ElSt98}, \cite[pp.~394--397]{PlPoStTa18}). Let
\mbox{$N \in \mathbb N$} with \mbox{$N \gg 1$} and \mbox{$M_1,\, M_2 \in 2 \mathbb N$} be given. By $\I_{M_1}$ we denote the index set \mbox{$\{-\frac{M_1}{2},\,1-\frac{M_1}{2},\,\ldots,\,\frac{M_1}{2}-1\}$}. We consider an \emph{exponential sum}
\mbox{$f:\,[-\frac{1}{2},\,\frac{1}{2}]\to \mathbb C$} of the form
\begin{equation}
\label{eq:f(x)}
f(x) \coloneqq \sum_{k\in \I_{M_1}} f_k\, {\mathrm e}^{-2\pi {\mathrm i}Nv_k x}\,, \quad x \in \big[-\tfrac{1}{2},\,\tfrac{1}{2}\big]\,,
\end{equation}
where \mbox{$f_k \in \mathbb C$} are given coefficients and \mbox{$v_k \in \big[-\frac{1}{2},\,\frac{1}{2}\big]$}, \mbox{$k\in \I_{M_1}$}, are arbitrary nodes in the frequency domain. The parameter \mbox{$N \in \mathbb N$} is called \emph{nonharmonic
bandwidth} of the exponential sum \eqref{eq:f(x)}.

We assume that a linear combination \eqref{eq:f(x)} of exponentials with bounded frequencies is given.
For arbitrary nodes \mbox{$x_j \in \big[-\frac{1}{2},\,\frac{1}{2}\big]$}, \mbox{$j \in \I_{M_2}$}, in the space
domain, we are interested in a fast evaluation of the $M_2$ values
\begin{equation}
\label{eq:f(xj)}
f(x_j) = \sum_{k\in \I_{M_1}} f_k\, {\mathrm e}^{-2\pi {\mathrm i}Nv_k x_j}\,, \quad j \in \I_{M_2}\,.
\end{equation}
A fast algorithm for the computation of the $M_2$ values \eqref{eq:f(xj)} is called a \emph{nonuniform fast Fourier transform with nonequispaced spatial and frequency data} (NNFFT) which was introduced by B.~Elbel and G.~Steidl in \cite{ElSt98}.
In this approach, the rapid evaluation of NNFFT is mainly based on the use of two compactly supported, continuous window functions.
As in \cite{LeGr05} this approach is also referred to as NFFT \emph{of type} 3.

In this paper we present new error estimates for the NNFFT.
Since these estimates depend exclusively on the so-called window parameters of the NNFFT, this gives rise to an appropriate parameter choice.
The outline of this paper is as follows. In Section~\ref{Sec:NNFFT}, we introduce the special set $\Omega$ of continuous, even functions \mbox{$\omega:\,\mathbb R \to [0,\,1]$} with the support \mbox{$[-1,\,1]$}.
Choosing \mbox{$\omega_1$, $\omega_2 \in \Omega$}, we consider two window functions
$$
\varphi_1(t) = \omega_1\bigg(\frac{N_1 t}{m_1}\bigg)\,, \quad \varphi_2(t) = \omega_2\bigg(\frac{N_2 t}{m_2}\bigg)\,, \quad t \in \mathbb R\,,
$$
where \mbox{$N_1 \coloneqq \sigma_1 N \in 2 \mathbb N$} with some oversampling factor \mbox{$\sigma_1>1$} and where \mbox{$m_1 \in {\mathbb N}\setminus \{1\}$} is a truncation parameter with \mbox{$2 m_1 \ll N_1$}.
Analogously, \mbox{$N_2 \coloneqq \sigma_2\,(N_1 + 2m_1) \in 2 \mathbb N$} is given with some oversampling factor \mbox{$\sigma_2 > 1$} and \mbox{$m_2 \in {\mathbb N}\setminus \{1\}$} is another truncation parameter with \mbox{$2 m_2 \ll \big(1 - \frac{1}{\sigma_1}\big)\,
N_2$}. For the fast, approximate computation of the values \eqref{eq:f(xj)}, we formulate the NNFFT in Algorithm~\ref{Alg2.1}. In Section~\ref{Sec:errorNNFFT}, we derive new explicit error estimates of
the NNFFT with two general window functions $\varphi_1$ and $\varphi_2$.
In Section~\ref{Sec:NNFFTsinhwindow}, we specify the result when using two $\sinh$-type window functions.
Namely, we show that for fixed nonharmonic bandwidth $N$ of \eqref{eq:f(x)}, the error of the related NNFFT has an exponential decay with respect to the truncation parameters $m_1$ and $m_2$. Numerical experiments illustrate the performance of our error estimates.

In Section~\ref{Sec:Approxsinc}, we study the approximation of the function  \mbox{$\mathrm{sinc}(N \pi x)$}, \mbox{$x \in [-1,\,1]$}, by an exponential sum. For given target accuracy \mbox{$\varepsilon > 0$} and \mbox{$n \ge 4 N$},
there exist coefficients \mbox{$w_j > 0$} and frequencies \mbox{$v_j \in (-1,\,1)$}, \mbox{$j =1\, \ldots,n$}, such that for all \mbox{$x \in [-1,\, 1]$},
\begin{equation*}
	\bigg|\mathrm{sinc}(N \pi x) - \sum_{j=1}^n w_j\,{\mathrm e}^{- \pi {\mathrm i}N v_j x}\bigg| \le \varepsilon\,.
\end{equation*}
In practice, we simplify the approximation procedure. 
Since for fixed \mbox{$N \in \mathbb N$}, it holds
$$
\mathrm{sinc} (N \pi x) = \frac{1}{2}\, \int_{-1}^1 {\mathrm e}^{- \pi {\mathrm i} N t x}\,{\mathrm d}t\,, \quad x \in \mathbb R\,,
$$
we apply the Clenshaw--Curtis quadrature with Chebyshev points \mbox{$z_k = \cos \frac{k \pi}{n}\in [-1,\,1]$}, \mbox{$k=0\, \ldots,n$}, where \mbox{$n\in \mathbb N$} fulfills \mbox{$n \ge 4 N$}.
Then the function \mbox{${\mathrm{sinc}}(N \pi x)$}, \mbox{$x \in [-1,\,1]$}, can be approximated by the exponential sum
\begin{equation}
\label{eq:sinc=expsum}
\sum_{k=0}^n w_k\,{\mathrm e}^{- \pi {\mathrm i} N z_k x}
\end{equation}
with explicitly known coefficients \mbox{$w_k > 0$} which satisfy the condition \mbox{$\sum_{k=0}^n w_k = 1$}.

An interesting signal processing application of the NNFFT is presented in the last Section~\ref{Sec:discreteSinc}.
If a signal \mbox{$h:\,\big[- \frac{1}{2},\, \frac{1}{2}\big] \to \mathbb C$} is to be reconstructed from its nonuniform
samples at \mbox{$a_k\in \big[- \frac{1}{2},\, \frac{1}{2}\big]$}, then $h$ is often modeled as linear combination of shifted $\mathrm{sinc}$ functions
$$
h(x) = \sum_{k \in \I_{L_1}} c_k \, \mathrm{sinc}\big(N \pi\,(x-a_k)\big)
$$
with complex coefficients $c_k$.
Hence, we present a fast, approximate computation of the \emph{discrete} $\mathrm{sinc}$ \emph{transform} (see \cite{GrLeIn06, LiBr11})
\begin{equation*}
h(b_{\ell}) = \sum_{k \in \I_{L_1}} c_k \, \mathrm{sinc}\big(N \pi\,(b_{\ell} - a_k)\big)\,, \quad \ell \in \I_{L_2}\,,
\end{equation*}
where \mbox{$b_{\ell} \in \big[- \frac{1}{2},\, \frac{1}{2}\big]$} can be nonequispaced.
The discrete $\mathrm{sinc}$ transform is motivated by numerous applications in signal processing.
However, since the $\mathrm{sinc}$ function decays slowly, it is often avoided in favor of some more local approximation. Here we prefer the approximation of the $\mathrm{sinc}$ function by an exponential
sum \eqref{eq:sinc=expsum}.
Then we obtain the fast $\mathrm{sinc}$ transform in Algorithm~\ref{alg:fastsinc}, which is an approximate algorithm for the fast computation of the values \eqref{eq:discretesinctransf} and applies the NNFFT twice.
Besides, the error of the fast $\mathrm{sinc}$ transform is estimated and numerical examples are presented as well.

\section{NNFFT}
\label{Sec:NNFFT}

\new{Now we start with the explanation of the main algorithm, the NNFFT.
To this end, we firstly introduce the special set $\Omega$, which is necessary to define required window functions $\varphi_j$, \mbox{$j=1,2$}.
Since the NNFFT is mainly based on the \mbox{well-known} NFFT, then we proceed with a short description of the NFFT and move on to the NNFFT afterwards.
This procedure is summarized in Algorithm~\ref{Alg2.1}.
Note that here a parameter \mbox{$a>1$} is necessary in order to prevent aliasing artifacts, since we approximate a non-periodic function on the interval \mbox{$[-1,\,1]$} by means of $a$-periodic functions.}

Let $\Omega$ be the set of all functions \mbox{$\omega:\, \mathbb R \to [0,\,1]$} with the following properties:
\begin{quote}
$\bullet$ Each function $\omega$ is even, has the support \mbox{$[-1,\,1]$}, and is continuous on $\mathbb R$.\\
$\bullet$ Each restricted function \mbox{$\omega|_{ [0,\, 1]}$} is decreasing with \mbox{$\omega(0) = 1$}.\\
$\bullet$ For each function $\omega$ its Fourier transform
$$
{\hat \omega}(v) \coloneqq \int_{\mathbb R} \omega(x)\,{\mathrm e}^{- 2 \pi {\mathrm i}\,v x}\,{\mathrm d}x = 2\,\int_0^1 \omega(x)\, \cos(2\pi v x)\,{\mathrm d}x
$$
is positive and decreasing for all \mbox{$v\in \big[0,\,\frac{m_1}{2 \sigma_1}\big]$}, where it holds \mbox{$m_1 \in \mathbb N \setminus \{1\}$} and \mbox{$\sigma_1 \in \big[\frac{5}{4}, \,2\big]$}.
\end{quote}
Obviously, each \mbox{$\omega \in \Omega$} is of bounded variation over \mbox{$[- 1,\, 1]$}.

\begin{example}
\label{Example:omega}
By $B_{2m_1}$, we denote the centered cardinal B-spline of even order $2m_1$ with \mbox{$m_1 \in \mathbb N$}. Thus, $B_2$ is the centered hat function. We consider the spline
$$
\omega_{\mathrm{B},1}(x) \coloneqq \frac{1}{B_{2m_1}(0)}\, B_{2m_1}(m_1 x)\,, \quad x \in \mathbb R\,,
$$
which has the support \mbox{$[-1,\,1]$}. Its Fourier transform reads as
$$
{\hat \omega}_{\mathrm{B},1}(v) = \frac{1}{m_1\, B_{2m_1}(0)}\,\Big(\mathrm{sinc}\,\frac{\pi v}{m_1}\Big)^{2m_1}\,, \quad v \in \mathbb R\,.
$$
Obviously, \mbox{${\hat \omega}_{\mathrm{B},1}(v)$} is positive and decreasing for \mbox{$v \in [0,\,m_1)$}. Hence, the function $\omega_{\mathrm{B},1}$ belongs to the set $\Omega$.

For \mbox{$\sigma_1 > \frac{\pi}{3}$} and \mbox{$\beta_1 = 3 m_1$} with \mbox{$m_1 \in \mathbb N \setminus \{1\}$}, we consider
$$
\omega_{\mathrm{alg},1}(x) \coloneqq \left\{ \begin{array}{ll} (1 - x^2)^{\beta_1 - 1/2} & \quad x \in [-1,\,1]\,,\\
0 & \quad x \in \mathbb R \setminus [-1,\,1]\,.
\end{array} \right.
$$
By $\mathrm{\cite[p.~8]{Oberh90}}$, its Fourier transform reads as
$$
{\hat \omega}_{\mathrm{alg},1}(v) = \frac{\pi\,(2 \beta_1)!}{4^{\beta_1}\,\beta_1!} \cdot \left\{ \begin{array}{ll}
(\pi v)^{-\beta_1}\, J_{\beta_1}(2 \pi v) & \quad v \in \mathbb R \setminus \{0\}\,, \\ [1ex]
\frac{1}{\beta_1!} & \quad v = 0\,,
\end{array} \right.
$$
where $J_{\beta_1}$ denotes the \emph{Bessel function of order} $\beta_1$. By $\mathrm{\cite[p.~370]{abst}}$, it holds for \mbox{$v \neq 0$} the equality
$$
(\pi v)^{-\beta_1}\, J_{\beta_1}(2 \pi v) = \frac{1}{\beta_1 !}\, \prod_{s=1}^{\infty} \bigg( 1 - \frac{4 \pi^2 v^2}{j_{\beta_1,s}^2}\bigg)\,,
$$
where $j_{\beta_1,s}$ denotes the $s$th positive zero of $J_{\beta_1}$.
For \mbox{$\beta_1 = 3 m_1$}, it holds \mbox{$j_{\beta_1,1} > 3 m_1 + \pi - \frac{1}{2}$} (see $\mathrm{\cite{IfSi85}}$).
Hence, by \mbox{$\sigma_1 > \frac{\pi}{3}$} we get
$$
\frac{2 \pi m_1}{2 \sigma_1 \,j_{\beta_1,1}} < \frac{\frac{\pi}{\sigma_1}\,m_1}{3 m_1 + \pi - \frac{1}{2}} <  \frac{3 m_1}{3 m_1 + \pi - \frac{1}{2}} < 1\,.
$$
Therefore, the Fourier transform ${\hat \omega}_{\mathrm{alg},1}(v)$ is positive and decreasing for \mbox{$v \in \big[0,\, \frac{m_1}{2 \sigma_1}\big]$}.
Hence, $\omega_{\mathrm{alg},1}$ belongs to the set $\Omega$.

Let \mbox{$\sigma_1 \in \big[\frac{5}{4},\,2 \big]$} and \mbox{$m_1 \in \mathbb N \setminus \{1\}$} be given. We consider the function
\begin{equation*}
\omega_{\sinh,1} (x) \coloneqq \left\{ \begin{array}{ll}
\frac{1}{\sinh \beta_1}\,\sinh \big(\beta_1 \sqrt{1 - x^2}\big)   &\quad  x \in [ - 1,\, 1]\,,\\
0                   & \quad x \in \mathbb R \setminus [ - 1,\, 1]
\end{array} \right.
\end{equation*}
with the \emph{shape parameter}
$$
\beta_1 \coloneqq 2 \pi m_1 \Big(1 - \frac{1}{2\sigma_1}\Big)\,.
$$
Then by $\mathrm{\cite[p.~38]{Oberh90}}$, its Fourier transform reads as
\begin{equation}
\label{eq:hatomegasinh}
{\hat \omega}_{\sinh,1}(v)
= \frac{\pi \beta_1}{\sinh \beta_1}\cdot \left\{ \begin{array}{ll}
(\beta_1^2 - 4\pi^2 v^2)^{-1/2}\,I_1\big(\sqrt{\beta_1^2 - 4 \pi^2 v^2}\big) &\quad |v| < m_1 \big(1 - \frac{1}{2\sigma_1}\big)\,,\\ [1ex]
\frac{1}{2} & \quad v = \pm m_1 \big(1 - \frac{1}{2\sigma_1}\big)\,, \\ [1ex]
(4 \pi^2 v^2 -\beta_1^2)^{-1/2}\,J_1\big(\sqrt{4 \pi^2 v^2 - \beta_1^2}\big) &\quad |v| > m_1 \big(1 - \frac{1}{2\sigma_1}\big)\,,
\end{array} \right.
\end{equation}
where $I_1$ and $J_1$ denote the \emph{modified Bessel function} and the \emph{Bessel function of first order}, respectively. Using the power series expansion of $I_1$ (see $\mathrm{\cite[p.~375]{abst}}$), we obtain for
\mbox{$|v| < m_1 \big(1 - \frac{1}{2\sigma_1}\big)$} that
$$
(\beta_1^2 - 4 \pi^2 v^2)^{-1/2}\,I_1\Big(\sqrt{\beta_1^2 - 4 \pi^2 v^2}\,\Big) = \frac{1}{2}\,\sum_{k=0}^{\infty} \frac{1}{4^k k! (k+1)!} \,(\beta_1^2 - 4 \pi^2 v^2)^k \,.
$$
Therefore, the Fourier transform ${\hat \omega}_{\sinh,1}(v)$ is positive and decreasing for \mbox{$v \in \big[0,\,\frac{m_1}{2 \sigma_1}\big]$}, since for \mbox{$\sigma_1 \ge \frac{5}{4}$} it holds
$$
\frac{m_1}{2 \sigma_1} < m_1 \Big(1 - \frac{1}{2\sigma_1}\Big)\,.
$$
Hence, $\omega_{\sinh,1}$ belongs to the set $\Omega$. $\Box$
\end{example}
\medskip

As known (see \cite{ElSt98, postta01}), the NNFFT can mainly be computed by means of an NFFT.
\new{This is why this algorithm is briefly explained below.}
For fixed \mbox{$N,\,\new{M_2} \in 2 \mathbb N$} and \mbox{$N_1 \coloneqq \sigma_1 N$} with \mbox{$\sigma_1>1$}, the NFFT (see \cite{DuRo93, duro95, St97} or \cite[pp.~377--381]{PlPoStTa18}) is a fast algorithm that approximately computes the values $p(x_j)$, \mbox{$j\in \I_{\new{M_2}}$}, of any 1-periodic trigonometric polynomial
\begin{align}
\label{eq:trigpoly}
	p(x) \coloneqq \sum_{k\in \I_N} c_k\,{\mathrm e}^{2\pi {\mathrm i}\,k x}
\end{align}
at nonequispaced nodes \mbox{$x_j \in \big[- \frac{1}{2}, \,\frac{1}{2}\big]$}, \mbox{$j \in \I_{\new{M_2}}$}, where \mbox{$c_k \in \mathbb C$, $k\in \I_N$}, are given complex coefficients.
\new{In other words, for the NFFT it holds \mbox{$N=M_1\in 2\mathbb N$} in \eqref{eq:f(xj)}.}

For \mbox{$\omega_1 \in \Omega$}
we introduce the \emph{window function}
\begin{equation}
\label{eq:firstwindow}
\varphi_1(t) \coloneqq \omega_1\bigg(\frac{N_1 t}{m_1}\bigg)\,, \quad t\in \mathbb R\,.
\end{equation}
By construction, the window function \eqref{eq:firstwindow} is even, has the support \mbox{$\big[- \frac{m_1}{N_1},\, \frac{m_1}{N_1}\big]$}, and is continuous on $\mathbb R$.
Further, the restricted window function \mbox{$\varphi_1 |_{[0,\,m_1/N_1]}$} is decreasing with \mbox{$\varphi_1(0) = 1$.}
Its Fourier transform
$$
{\hat \varphi}_1(v) \coloneqq \int_{\mathbb R} \varphi_1(t)\,{\mathrm e}^{- 2 \pi {\mathrm i}\,v t}\,{\mathrm d}t = 2\,\int_0^{m_1/N_1} \varphi_1(t)\, \cos(2\pi v t)\,{\mathrm d}t
$$
is positive and decreasing for \mbox{$v \in \big[0,\,N_1 - \frac{N}{2}\big)$}.
Thus, $\varphi_1$ is of bounded variation over \mbox{$\big[-\frac{1}{2}, \, \frac{1}{2}\big]$}.

In the following, we denote the torus $\mathbb R/\mathbb Z$ by $\mathbb T$ and the Banach space  of continuous, \mbox{1-periodic} functions by $C(\mathbb T)$.
For the window function \eqref{eq:firstwindow}, we denote its \mbox{1-periodization} by
$$
{\tilde \varphi}_1^{(1)}(x) \coloneqq \sum_{k \in \mathbb Z} \varphi_1(x + k)\,, \quad x \in \mathbb R\,.
$$

Using a linear combination of shifted versions of the \mbox{1-periodized} window function ${\tilde \varphi}_1^{(1)}$, we construct a \mbox{1-periodic} continuous function \mbox{$s \in C(\mathbb T)$} which approximates \eqref{eq:trigpoly} well.
Then the computation of the values $s(x_j)$, \mbox{$j \in \I_{\new{M_2}}$}, is very easy, since $\varphi_1$ has the small support \mbox{$\big[-\frac{m_1}{N_1},\, \frac{m_1}{N_1}\big]$}.
The computational cost of NFFT is \mbox{${\mathcal O}\big(N \log N + \new{M_2}\big)$} flops, see \cite{DuRo93, duro95, St97} or \cite[pp.~377--381]{PlPoStTa18}.
The error of the NFFT (see \cite{PoTa20}) can be estimated by
\begin{eqnarray*}
\max_{j \in \I_{\new{M_2}}} |s(x_j) - p(x_j)| &\le& \|s - p\|_{C(\mathbb T)} \coloneqq \max_{x\in [-1/2,1/2]} |s(x) - p(x)|\\
&\le& e_{\sigma_1}(\varphi_1)\,\sum_{n \in \I_N} |c_n|\,,
\end{eqnarray*}
where $e_{\sigma_1}(\varphi_1)$ denotes the $C(\mathbb T)$-\emph{error constant} defined as
\begin{equation}
\label{eq:esigma(varphi)}
e_{\sigma_1}(\varphi_1) = \sup_{N \in 2 \mathbb N} e_{\sigma_1,N}(\varphi_1)
\end{equation}
with
$$
e_{\sigma_1,N}(\varphi_1) \coloneqq \max_{n\in \I_N} \Bigg\| \sum_{r \in \mathbb Z \setminus \{0\}} \frac{{\hat \varphi}_1(n + r N_1)}{{\hat \varphi}_1(n)}\,{\mathrm e}^{2 \pi {\mathrm i}\,rN_1 \cdot\,}\Bigg\|_{C(\mathbb T)}\,.
$$
Note that the constants $e_{\sigma_1,N}(\varphi_1)$ are bounded with respect to $N$ (see \cite[Theorem~ 5.1]{PoTa20}).
\medskip

\new{Now we proceed with the NNFFT.
For better readability, we describe the procedure just shortly. For more detailed explanations we refer to \cite{ElSt98}.}
For chosen functions \mbox{$\omega_1,\, \omega_2 \in \Omega$}, we form the window functions
\begin{equation}
\label{eq:twowindows}
\varphi_1(t) \coloneqq \omega_1\bigg(\frac{N_1 t}{m_1}\bigg)\,, \quad \varphi_2(t) \coloneqq \omega_2\bigg(\frac{N_2 t}{m_2}\bigg)\,, \quad t\in \mathbb R\,,
\end{equation}
where again \mbox{$N_1 \coloneqq \sigma_1 N \in 2 \mathbb N$} with some oversampling factor \mbox{$\sigma_1 > 1$} and \mbox{$m_1 \in \mathbb N \setminus \{1\}$} with \mbox{$2 m_1 \ll N_1$} and where
\mbox{$N_2\coloneqq \sigma_2\,(N_1 + 2 m_1) \in 2 \mathbb N$}
with an oversampling factor \mbox{$\sigma_2 > 1$} and \mbox{$m_2 \in \mathbb N \setminus \{1\}$} with \mbox{$2 m_2 \le \big(1 - \frac{1}{\sigma_1}\big)\, N_2$}.
The second window function $\varphi_2$ has the support \mbox{$\big[-\frac{m_2}{N_2},\,\frac{m_2}{N_2}\big]$}.
\new{Additionally, in order to prevent aliasing, we use $a$-periodic functions, where} we introduce the constant
\begin{equation}
\label{eq:consta}
a\coloneqq 1 + \frac{2m_1}{N_1} > 1\,,
\end{equation}
such that \mbox{$a N_1 = N_1 + 2m_1$} and \mbox{$N_2 = \sigma_2 \sigma_1 a N$}.
Without loss of generality, we can assume that
\begin{equation}
\label{eq:condvk}
v_k \in \big[- \tfrac{1}{2a}, \, \tfrac{1}{2a}\big]\,.
\end{equation}
If \mbox{$v_k \in \big[-\frac{1}{2},\,\frac{1}{2}\big]$},
then we replace the nonharmonic bandwidth $N$ by \mbox{$N^{\ast} \coloneqq N + \lceil \frac{2m_1}{\sigma_1}\rceil$} and set \mbox{$v_j^{\ast} \coloneqq \frac{N}{N^{\ast}}\,v_j \in \big[- \frac{1}{2a}, \, \frac{1}{2a}\big]$} such that \mbox{$N v_j = N^{\ast} v_j^{\ast}$}.

For arbitrarily given \mbox{$f_k \in \mathbb C$}, \mbox{$k \in \I_{M_1}$}, and
\mbox{$v_k \in \big[-\frac{1}{2a},\,\frac{1}{2a}\big]$}, \mbox{$k\in \I_{M_1}$},
we introduce the compactly supported, continuous auxiliary function
$$
h(t) \coloneqq \sum_{k \in \I_{M_1}} f_k\,\varphi_1(t - v_k)\,, \quad t \in \mathbb R\,,
$$
which has the Fourier transform
\begin{eqnarray}
\hat h (N x) &=& \int_{\mathbb R} h(t)\,{\mathrm e}^{-2 \pi {\mathrm i}N  x t}\,{\mathrm d}t \nonumber \\
&=& \sum_{k \in \I_{M_1}} f_k \, \int_{\mathbb R} \varphi_1(t - v_k)\,{\mathrm e}^{-2 \pi {\mathrm i}N x t}\,{\mathrm d}t \label{eq:hath(v)conv} \\
&=& \sum_{k \in \I_{M_1}} f_k\,{\mathrm e}^{-2 \pi {\mathrm i}N v_k x}\,{\hat \varphi}_1(N x) = f(x)\,{\hat \varphi}_1(N x)\,, \quad x \in \mathbb R\,. \label{eq:hath(v)}
\end{eqnarray}
Hence, for arbitrary nodes $x_j \in \big[- \frac{1}{2}, \, \frac{1}{2}\big]$, $j\in \I_{M_2}$, we have
$$
f(x_j) = \frac{{\hat h}(Nx_j)}{{\hat \varphi}_1(Nx_j)}\,, \quad j \in \I_{M_2}\,.
$$
Therefore, it remains to compute the values ${\hat h}(Nx_j)$, \mbox{$j\in \I_{M_2}$}, because we can precompute the values ${\hat \varphi}_1(Nx_j)$, \mbox{$j \in \I_{M_2}$}.
In some cases (see Section~4), these values ${\hat \varphi}_1(Nx_j)$, \mbox{$j\in \I_{M_2}$}, are explicitly known.

For arbitrary \mbox{$v_k \in \big[- \frac{1}{2a},\,\frac{1}{2a}\big]$}, \mbox{$k \in \I_{M_1}$},
we have \mbox{$\varphi_1(t - v_k) = 0$} for all
\mbox{$t < -\frac{1}{2a} - \frac{m_1}{N_1} =$} \mbox{$-\frac{a}{2} +\big(\frac{1}{2}- \frac{1}{2a}\big)$}
and for all
\mbox{$t > \frac{1}{2a} + \frac{m_1}{N_1} = \frac{a}{2} - \big(\frac{1}{2}- \frac{1}{2a}\big)$},
since \mbox{${\mathrm{supp}}\, \varphi_1 = \big[- \frac{m_1}{N_1},\,\frac{m_1}{N_1}\big]$}
and \mbox{$\frac{1}{2}- \frac{1}{2a}> 0$}.
Thus, by \eqref{eq:hath(v)conv} and
$$
{\mathrm{supp}}\,\varphi_1(\cdot - v_k) \subset \bigg[-\frac{a}{2},\, \frac{a}{2}\bigg]\,, \quad k\in \I_{M_1}\,,
$$
we obtain
$$
\hat h (N x) = \sum_{k\in \I_{M_1}} f_k\, \int_{-a/2}^{a/2} \varphi_1(t-v_k)\,{\mathrm e}^{-2 \pi {\mathrm i}N x t}\,{\mathrm d}t\,, \quad x \in \mathbb R\,.
$$
Then the rectangular quadrature rule leads to
\begin{equation}
\label{eq:s(Nx):=}
s(N x) \coloneqq \sum_{k\in \I_{M_1}} f_k\,\frac{1}{N_1} \sum_{\ell \in \I_{N_1 + 2m_1}} \varphi_1\Big(\frac{\ell}{N_1} - v_k\Big)\, {\mathrm e}^{-2 \pi {\mathrm i}\,\ell x/\sigma_1}\,, \quad x \in \mathbb R\,,
\end{equation}
which approximates ${\hat h}(N x)$.
Note that $\frac{\ell}{N_1} \in \big[- \frac{a}{2},\,\frac{a}{2}\big]$ for each $\ell \in \I_{N_1 + 2 m_1}$ by \mbox{$N_1 + 2m_1 = a N_1$}.
Changing the order of summations in \eqref{eq:s(Nx):=}, it follows that
\begin{equation}
\label{eq:s(v):=}
s(N x) = \sum_{\ell \in \I_{N_1 + 2m_1}} \Bigg(\frac{1}{N_1} \sum_{k\in \I_{M_1}} f_k\, \varphi_1\Big(\frac{\ell}{N_1} - v_k\Big)\Bigg)\, {\mathrm e}^{-2 \pi {\mathrm i}\,\ell x/\sigma_1}\,.
\end{equation}
After computation of the inner sums
\begin{equation}
\label{eq:gell}
g_{\ell} \coloneqq \frac{1}{N_1} \sum_{k\in \I_{M_1}} f_k\, \varphi_1\Big(\frac{\ell}{N_1} - v_k\Big)\,, \quad \ell \in \I_{N_1 + 2m_1}\,,
\end{equation}
we arrive at the following NFFT
$$
s(N \,x_j) = \sum_{\ell \in \I_{N_1 + 2m_1}} g_{\ell}\, {\mathrm e}^{-2 \pi {\mathrm i}\,\ell x_j/\sigma_1}\,, \quad j \in \I_{M_2}\,.
$$
If we denote the result of this NFFT  (with the 1-periodization ${\tilde \varphi}_2^{(1)}$ of the second window function $\varphi_2$ and \mbox{$N_2 \coloneqq \sigma_2\,(N_1 + 2 m_1)$}) by $s_1(Nx_j)$, then \mbox{$s_1(Nx_j)/{\hat \varphi}_1(Nx_j)$} is an approximate value of $f(x_j)$, \mbox{$j \in \I_{M_2}$}.
\new{Thus, the algorithm can be summarized as follows.}

\begin{algorithm}[NNFFT]\phantom{.}
\label{Alg2.1}

\textit{Input}: \new{Nonharmonic bandwidth} \mbox{$N \in \mathbb N$} with \mbox{$N \gg 1$}, \new{numbers of nodes} \mbox{$M_1,\, M_2\in 2{\mathbb N}$},
\mbox{$N_1 \coloneqq \sigma_1N\in 2{\mathbb N}$} with \new{oversampling factor} \mbox{$\sigma_1> 1$} and \new{truncation parameter} \mbox{$m_1\in {\mathbb N}\setminus \{1\}$} with \mbox{$2m_1 \ll N_1$},
\mbox{$N_2 \coloneqq \sigma_2\,(N_1 + 2m_1)\in 2{\mathbb N}$} with \new{oversampling factor} \mbox{$\sigma_2 > 1$} and \new{truncation parameter} \mbox{$m_2\in {\mathbb N}\setminus \{1\}$} with \mbox{$2m_2 \le \big(1 - \frac{1}{\sigma_1}\big)\, N_2$},\\
\new{arbitrary nodes} \mbox{$v_k \in \big[-\frac{1}{2a},\,\frac{1}{2a}\big]$}, \mbox{$k\in \I_{M_1}$}, \new{in the frequency domain} with \mbox{$a \coloneqq 1 + \frac{2m_1}{N_1}$},
\new{arbitrary nodes} \mbox{$x_j \in \big[-\frac{1}{2},\,\frac{1}{2}\big]$}, \mbox{$j \in \I_{M_2}$}, \new{in the spatial domain}
\new{as well as window functions} $\varphi_1$ and $\varphi_2$ given by \eqref{eq:twowindows}.
\medskip

$0$. Precompute the following values: \vspace*{-1.2ex}
\begin{itemize}
	\item[\new{(i)}] ${\hat \varphi}_1(Nx_j)$ for $j\in \I_{M_2}$, ${\hat \varphi}_2\big(\frac{\ell}{a}\big)$ for $\ell \in \I_{N_1+2m_1}$, \vspace*{-1.2ex}
	\item[\new{(ii)}] $\varphi_1\big(\frac{\ell}{N_1} - v_k\big)$ for $k\in \I_{M_1}$ and \mbox{$\ell \in \I_{N_1+2m_1}'(v_k)\coloneqq\{\ell \in \I_{N_1+2m_1}:\, \big|\frac{\ell}{N_1} - v_k\big| < \frac{m_1}{N_1}\}$}, \vspace*{-3.7ex}
	\item[\new{(iii)}] $\varphi_2\big(\frac{x_j}{\sigma_1} - \frac{s}{N_2}\big)$ for $j \in \I_{M_2}$ and \mbox{$s\in \I_{N_2}''(x_j)\coloneqq\{s\in \I_{N_2}:\, \big|\frac{s}{N_2} - \frac{x_j}{\sigma_1}\big| < \frac{m_2}{N_2}\}$}, \vspace*{-1.2ex}
	\item[\new{(iv)}] Further set  $\varphi_1\big(\frac{\ell}{N_1} - v_k\big)\coloneqq0$ for $k\in \I_{M_1}$ and \mbox{$\ell \in \I_{N_1+2m_1}\setminus \I_{N_1+2m_1}'(v_k)$}.
\end{itemize}

$1$. For all $\ell \in \I_{N_1+2m_1}$ compute the sums \eqref{eq:gell}.\\
$2$. For all $\ell \in \I_{N_1 + 2m_1}$ form the values
$$
{\hat g}_{\ell} \coloneqq \frac{g_{\ell}}{{\hat \varphi}_2(\ell)}\,.
$$
$3$. For all $s\in \I_{N_2}$ compute by fast Fourier transform $\mathrm{(FFT)}$ of length $N_2$
$$
h_s \coloneqq \frac{1}{N_2}\,\sum_{\ell \in \I_{N_1+2m_1}} {\hat g}_{\ell}\,{\mathrm e}^{-2 \pi {\mathrm i}\,\ell s/N_2}\,.
$$
$4$. For all $j\in \I_{M_2}$ calculate the short sums
$$
s_1(N x_j) \coloneqq \sum_{s \in \I_{N_2}''(x_j)} h_s\,\varphi_2\Big(\frac{x_j}{\sigma_1} - \frac{s}{N_2}\Big)\,.
$$
\smallskip

{\textit Output}: $s_1(Nx_j)/{\hat \varphi}_1(Nx_j)$ approximate value of \eqref{eq:f(xj)} for $j\in \I_{M_2}$.
\end{algorithm}

The computational cost of the NNFFT is equal to \mbox{${\mathcal O}\big(N \log N + M_1 + M_2\big)$} flops.

In Step 4 of Algorithm~\ref{Alg2.1} we use the assumption $2m_2 \le \big(1 - \frac{1}{\sigma_1}\big)\, N_2$ such that
$$
\frac{1}{2 \sigma_1} + \frac{m_2}{N_2} \le \frac{1}{2}\,.
$$
Then for all $j \in \I_{M_2}$ and $s\in \I_{N_2}$, it holds
$$
{\tilde \varphi}_2^{(1)}\Big(\frac{x_j}{\sigma_1} - \frac{s}{N_2}\Big) = \varphi_2\Big(\frac{x_j}{\sigma_1} - \frac{s}{N_2}\Big)\,.
$$
\new{Since we approximate a non-periodic function $f$ on the interval \mbox{$\big[-\frac{1}{2},\,\frac{1}{2}\big]$} by means of \mbox{$a$-periodic} functions on the torus \mbox{$a\mathbb T\cong\big[-\frac{a}{2},\,\frac{a}{2}\big)$}, the parameter $a$ has to fulfill the condition \mbox{$a>1$}, in order to prevent aliasing artifacts.}

\section{Error estimates for NNFFT}
\label{Sec:errorNNFFT}

Now we study the error of the NNFFT, which is measured in the form
\begin{equation*}
	\max_{j \in \I_{M_2}} \bigg| f(x_j) - \frac{s_1(Nx_j)}{{\hat \varphi}_1(Nx_j)}\bigg|\,,
\end{equation*}
where $f$ is a given exponential sum \eqref{eq:f(x)} and $x_j \in \big[-\frac{1}{2},\,\frac{1}{2}\big]$, $j \in \I_{M_2}$, are arbitrary spatial nodes.
\new{At the beginning of this section we present some technical lemmas. The main result will be Theorem~\ref{Thm:errorestimate}.}

We introduce the $a$-{\emph{periodization}} of the given window function \eqref{eq:firstwindow} by
\begin{equation}
\label{eq:aperiodicwindow}
{\tilde \varphi}_1^{(a)}(x) \coloneqq \sum_{\ell \in \mathbb Z} \varphi_1(x + a\, \ell)\,, \quad x \in \mathbb R\,.
\end{equation}
For each $x\in \mathbb R$, the above series \eqref{eq:aperiodicwindow} has at most one nonzero term. This can be seen as follows: For arbitrary $x \in \mathbb R$ there exists a unique $\ell^* \in \mathbb Z$ such
that $x = - a\,\ell^* + r$ with a residuum $r \in \big[-\frac{a}{2},\,\frac{a}{2}\big)$. Then $\varphi_1(x + a\,\ell^*) = \varphi_1(r)$ and hence $\varphi_1(r) > 0$ for $r \in \big(-\frac{m_1}{N_1},\,\frac{m_1}{N_1}\big)$ and $\varphi_1(r) = 0$ for
$r \in \big[-\frac{a}{2},\,- \frac{m_1}{N_1}\big] \cup \big[ \frac{m_1}{N_1},\,\frac{a}{2}\big)$. For each $\ell \in \mathbb Z \setminus \{\ell^*\}$, we have
$$
\varphi_1 (x + a\, \ell) = \varphi_1\big(a \,(\ell - \ell^*) + r\big) = 0\,,
$$
since $\big|a \,(\ell - \ell^*) + r\big| \ge \frac{a}{2} = \frac{1}{2} + \frac{m_1}{N_1} > \frac{m_1}{N_1}$. Further it holds
$$
{\tilde \varphi}_1^{(a)}(x) = \varphi_1(x)\,, \quad x \in \big[-1-\tfrac{m_1}{N_1},\,1 + \tfrac{m_1}{N_1}\big]\,.
$$
By the construction of $\varphi_1$, the $a$-periodic window function \eqref{eq:aperiodicwindow} is continuous on $\mathbb R$ and of bounded variation over $\big[- \frac{a}{2},\, \frac{a}{2}\big]$.
Then the $k$th Fourier coefficient of the $a$-periodic window function \eqref{eq:aperiodicwindow} reads as follows
\begin{equation}
\label{eq:Fouriercoeff}
c_k^{(a)}\Big(\!{\tilde \varphi}_1^{(a)}\!\Big) \coloneqq \frac{1}{a} \int_{-a/2}^{a/2} {\tilde \varphi}_1^{(a)}(t)\,{\mathrm e}^{-2 \pi {\mathrm i}\,k t/a}\,{\mathrm d}t = \frac{1}{a}\,{\hat \varphi}_1\bigg(\!\frac{k}{a}\!\bigg)\,, \quad k \in \mathbb Z\,.
\end{equation}
By the convergence theorem of Dirichlet--Jordan (see \cite[Vol.~1, pp.~57--58]{Zygmund}), the $a$\hbox{-}periodic Fourier series of \eqref{eq:aperiodicwindow} converges uniformly on $\mathbb R$ and it holds
\begin{equation}
\label{eq:Fouriervarphia}
{\tilde \varphi}_1^{(a)}(x) = \sum_{k\in \mathbb Z} c_k^{(a)}\Big(\!{\tilde \varphi}_1^{(a)}\!\Big) \,{\mathrm e}^{2 \pi {\mathrm i}\,k x/a}
= \frac{1}{a}\, \sum_{k\in \mathbb Z} {\hat \varphi}_1\bigg(\!\frac{k}{a}\!\bigg)\,{\mathrm e}^{2 \pi {\mathrm i}\,k x/a}\,.
\end{equation}
\new{Then we have the following technical lemma.}

\begin{Lemma}
\label{lemma:Fourierseries}
Let the window function $\varphi_1$ be given by \eqref{eq:firstwindow}. Then for any $n \in \I_N$ with $N \in 2 \mathbb N$, the series
$$
\sum_{r\in \mathbb Z} c_{n + r\,(N_1 + 2m_1)}^{(a)}\Big(\!{\tilde \varphi}_1^{(a)}\!\Big) \,{\mathrm e}^{2 \pi {\mathrm i}\,(n + r\,(N_1+ 2m_1)) x/a}
$$
is uniformly convergent on $\mathbb R$ and has the sum
$$
\frac{1}{N_1 + 2m_1} \sum_{\ell \in \I_{N_1+2m_1}} {\mathrm e}^{-2 \pi {\mathrm i}\, n \ell/(N_1 + 2m_1)}\, {\tilde \varphi}_1^{(a)}\bigg(\!x + \frac{\ell}{N_1}\!\bigg)
$$
which coincides with the rectangular quadrature rule of the integral
$$
c_n^{(a)}\Big(\!{\tilde \varphi}_1^{(a)}(x + \cdot)\!\Big) = \frac{1}{a} \int_{-a/2}^{a/2} {\tilde \varphi}_1^{(a)}(x+ s)\,{\mathrm e}^{2 \pi {\mathrm i}\,n s/a} \, {\mathrm d}s = c_n^{(a)}\Big(\!{\tilde \varphi}_1^{(a)}\!\Big)\, {\mathrm e}^{2 \pi {\mathrm i}\,n x/a}\,.
$$
\end{Lemma}

\emph{Proof.} Using the uniformly convergent Fourier series \eqref{eq:Fouriervarphia}, we obtain for all $n\in \I_N$ that
$$
{\mathrm e}^{-2 \pi {\mathrm i}\,n x/a}\,{\tilde \varphi}_1^{(a)}(x) = \sum_{k\in \mathbb Z} c_k^{(a)}\Big(\!{\tilde \varphi}_1^{(a)}\!\Big) \,{\mathrm e}^{2 \pi {\mathrm i}\,(k-n) x/a}
=\sum_{q\in \mathbb Z} c_{n+q}^{(a)}\Big(\!{\tilde \varphi}_1^{(a)}\!\Big) \,{\mathrm e}^{2 \pi {\mathrm i}\,q x/a}\,.
$$
Replacing $x$ by $x + \frac{\ell}{N_1}$ with $\ell \in \I_{N_1+ 2m_1}$, we see that by $N_1 + 2m_1 = aN_1$,
$$
{\mathrm e}^{-2 \pi {\mathrm i}\,n (x + \ell/N_1)/a}\,{\tilde \varphi}_1^{(a)}\bigg(\!x + \frac{\ell}{N_1}\!\bigg) =\sum_{q\in \mathbb Z} c_{n+q}^{(a)}\Big(\!{\tilde \varphi}_1^{(a)}\!\Big) \,{\mathrm e}^{2 \pi {\mathrm i}\,q x/a}\,
{\mathrm e}^{2\pi {\mathrm i}\,q \ell/(N_1 + 2m_1)}\,.
$$
Summing the above formulas for all $\ell\in \I_{N_1+ 2m_1}$ and applying the known formula
$$
\sum_{\ell\in \I_{N_1+2m_1}}{\mathrm e}^{2\pi {\mathrm i}\,q \ell/(N_1 + 2m_1)} = \left\{ \begin{array}{ll}
N_1+ 2m_1 & \quad q \equiv 0\,\mathrm{mod}\,(N_1 + 2m_1)\,,\\
0 & \quad q \not\equiv 0\, \mathrm{mod}\,(N_1 + 2m_1)\,,
\end{array} \right.
$$
we conclude that
\begin{eqnarray*}
& &\sum_{\ell\in \I_{N_1+2m_1}} {\mathrm e}^{-2 \pi {\mathrm i}\,n (x + \ell/N_1)/a}\,{\tilde \varphi}_1^{(a)}\bigg(\!x + \frac{\ell}{N_1}\!\bigg)\\
& & =\, (N_1 + 2m_1) \sum_{r\in \mathbb Z}  c_{n + r\,(N_1 + 2m_1)}^{(a)}\Big(\!{\tilde \varphi}_1^{(a)}\!\Big) \,{\mathrm e}^{2 \pi {\mathrm i}\, r\,(N_1+ 2m_1) x/a}\,.
\end{eqnarray*}
Obviously,
$$
\frac{1}{N_1+2m_1} \sum_{\ell\in \I_{N_1+2m_1}} {\mathrm e}^{-2 \pi {\mathrm i}\,n (x + \ell/N_1)/a}\,{\tilde \varphi}_1^{(a)}\bigg(\!x + \frac{\ell}{N_1}\!\bigg)
$$
is the rectangular quadrature formula of the integral
$$
\frac{1}{a} \int_{-a/2}^{a/2} {\tilde \varphi}_1^{(a)}(x+ s)\,{\mathrm e}^{2 \pi {\mathrm i}\,n s/a} \, {\mathrm d}s
$$
with respect to the uniform grid $\big\{\frac{\ell}{N_1} :\,\ell \in \I_{N_1+2m_1}\big\}$ of the interval $\big[-\frac{a}{2},\, \frac{a}{2}\big]$.
This completes the proof. \qedsymbol
\medskip

By \eqref{eq:Fouriercoeff} we obtain that for $n \in \I_N$,
$$
\Bigg| \sum_{r \in \mathbb Z \setminus \{0\}} \frac{c_{n + r\,(N_1+2m_1)}^{(a)}({\tilde \varphi}_1^{(a)})}{c_n^{(a)}({\tilde \varphi}_1^{(a)})}\,{\mathrm e}^{2 \pi {\mathrm i}\,r\,(N_1+2m_1) x/a}\Bigg|
= \Bigg| \sum_{r \in \mathbb Z \setminus \{0\}} \frac{{\hat \varphi}_1(n/a + r N_1)}{{\hat \varphi}_1(n/a)}\,{\mathrm e}^{2 \pi {\mathrm i}\,rN_1 x/a}\Bigg|\,.
$$

Now we generalize the \new{technical} Lemma~\ref{lemma:Fourierseries}.

\begin{Lemma}
\label{lemma:Fourierseriesnew}
For arbitrary fixed $v \in \big[-\frac{N}{2},\,\frac{N}{2}\big]$, $N\in\mathbb N$, and given window function~\eqref{eq:firstwindow}, the function
\begin{equation}
\label{eq:defpsi}
\psi_1(x)\coloneqq \frac{1}{N_1}\,\sum_{\ell\in \mathbb Z} {\mathrm e}^{- 2\pi {\mathrm i}\,v \ell/(N_1 + 2m_1)}\, {\mathrm e}^{-2 \pi {\mathrm i}\,v x/a}\,\varphi_1\bigg(\!x +\frac{\ell}{N_1}\!\bigg)
\end{equation}
is $\frac{1}{N_1}$-periodic, continuous on $\mathbb R$, and of bounded variation over $\big[-\frac{1}{2},\, \frac{1}{2}\big]$. For each $x\in \mathbb R$, the corresponding $\frac{1}{N_1}$-periodic Fourier series converges
uniformly to $\psi_1(x)$, i.\,e.,
$$
\psi_1(x) = \sum_{r \in \mathbb Z} {\hat \varphi}_1\Big(\frac{v}{a} + r N_1\Big)\,{\mathrm e}^{2 \pi {\mathrm i}\,rN_1 x}\,.
$$
\end{Lemma}

\emph{Proof.} The definition \eqref{eq:defpsi} of the function $\psi_1$ is correct, since
$$
\psi_1(x) =  \frac{1}{N_1}\,\sum_{\ell\in {\mathbb Z}_{m_1,N_1}(x)} {\mathrm e}^{- 2\pi {\mathrm i}\,v \ell/(N_1 + 2m_1)}\, {\mathrm e}^{-2 \pi {\mathrm i}\,v x/a}\,\varphi_1\bigg(\!x +\frac{\ell}{N_1}\!\bigg)
$$
with the finite index set ${\mathbb Z}_{m_1,N_1}(x) = \{ \ell\in \mathbb Z: \, | N_1 x + \ell| < m_1 \}$.
If $x\in [-\frac{1}{2},\,\frac{1}{2}\big]$, we observe that ${\mathbb Z}_{m_1,N_1}(x) \subseteq \I_{N_1 + 2m_1}$ and therefore
$$
\psi_1(x) =  \frac{1}{N_1}\,\sum_{\ell\in \I_{N_1+2m_1}} {\mathrm e}^{- 2\pi {\mathrm i}\,v \ell/(N_1 + 2m_1)}\, {\mathrm e}^{-2 \pi {\mathrm i}\,v x/a}\,\varphi_1\bigg(\!x +\frac{\ell}{N_1}\!\bigg)\,.
$$
Simple calculation shows that for each $x \in \mathbb R$,
$$
\psi_1\bigg(\!x + \frac{1}{N_1}\!\bigg) = \frac{1}{N_1}\,\sum_{\ell\in \mathbb Z} {\mathrm e}^{- 2\pi {\mathrm i}\,v (\ell+1)/(N_1 + 2m_1)}\, {\mathrm e}^{-2 \pi {\mathrm i}\,v x/a}\,\varphi_1\bigg(\!x +\frac{\ell+1}{N_1}\!\bigg) = \psi_1(x)\,.
$$
By the construction of $\varphi_1$, the $\frac{1}{N_1}$-periodic function $\psi_1$ is  continuous on $\mathbb R$ and of bounded variation over $\big[-\frac{1}{2}, \, \frac{1}{2}\big]$. Thus, by the convergence theorem of Dirichlet--Jordan, the
Fourier series of $\psi_1$ converges uniformly on $\mathbb R$ to $\psi_1$. The $r$th Fourier coefficient of $\psi_1$ reads as follows
\begin{eqnarray*}
c_r^{(1/N_1)}(\psi_1) &=& N_1\,\int_0^{1/N_1} \psi_1(t)\,{\mathrm e}^{-2 \pi {\mathrm i}\,rN_1 t}\, {\mathrm d}t\\
&=& \sum_{\ell\in \mathbb Z} {\mathrm e}^{- 2\pi {\mathrm i}\,v \ell/(N_1 + 2m_1)}\, \int_0^{1/N_1} {\mathrm e}^{-2 \pi {\mathrm i}\,v t/a}\,\varphi_1\bigg(\!t +\frac{\ell}{N_1}\!\bigg)\,{\mathrm d}t\\
&=& \sum_{\ell\in \mathbb Z} \int_{\ell/N_1}^{(\ell+1)/N_1} \varphi_1(s)\,{\mathrm e}^{-2 \pi {\mathrm i}\,(v/a + rN_1)\,s }\,{\mathrm d}s = {\hat \varphi}_1\Big(\frac{v}{a} + r N_1\Big)\,, \quad r \in \mathbb Z\,.
\end{eqnarray*}
This completes the proof. \qedsymbol
\medskip

From Lemma~\ref{lemma:Fourierseriesnew} \new{leads immediately to the following technical result}.

\begin{corollary}
\label{Cor:Cor2.3}
Let the window function $\varphi_1$ be given by \eqref{eq:firstwindow}. For all $x \in \big[-\frac{1}{2},\,\frac{1}{2}\big]$ and $w \in \big[-\frac{N}{2a},\,\frac{N}{2a}\big]$ it holds then
\begin{eqnarray}
\label{eq:sumhatvarphi}
& &\sum_{r \in {\mathbb Z}\setminus \{0\}} \frac{{\hat \varphi}_1(w + r N_1)}{{\hat \varphi}_1(w)}\,{\mathrm e}^{2 \pi {\mathrm i}\,(w + r N_1)\,x} \nonumber\\
& & =\,\frac{1}{N_1\,{\hat \varphi}_1(w)} \sum_{\ell\in \I_{N_1 + 2m_1}} {\mathrm e}^{- 2\pi {\mathrm i}\,w \ell/N_1}\,\varphi_1\bigg(\!x +\frac{\ell}{N_1}\!\bigg) - {\mathrm e}^{2 \pi {\mathrm i}\,w x}\,.
\end{eqnarray}
Further, for all $w \in \big[-\frac{N}{2a},\,\frac{N}{2a}\big]$, it holds
\begin{eqnarray}
\label{eq:maxmax}
& &\max_{x\in [-1/2,\,1/2]} \Bigg| \frac{1}{N_1\,{\hat \varphi}_1(w)}\sum_{\ell \in \I_{N_1+2m_1}} \varphi_1\bigg(\!x +\frac{\ell}{N_1}\!\bigg)\,
{\mathrm e}^{- 2\pi {\mathrm i}\,w \ell /N_1} - {\mathrm e}^{2 \pi {\mathrm i}\,w x}\Bigg| \nonumber\\
& &=\, \Bigg\|\sum_{r \in {\mathbb Z}\setminus \{0\}} \frac{{\hat \varphi}_1(w + r N_1)}{{\hat \varphi}_1(w)}\,{\mathrm e}^{2 \pi {\mathrm i}\, r N_1\,\cdot}\Bigg\|_{C(\mathbb T)}\,.
\end{eqnarray}
\end{corollary}

\emph{Proof.} As before, let $v \in \big[-\frac{N}{2},\,\frac{N}{2}\big]$ be given.
Substituting $w \coloneqq \frac{v}{a}\in \big[-\frac{N}{2a},\,\frac{N}{2a}\big]$ and observing $N_1 + 2 m_1 = a N_1$, we obtain by Lemma~\ref{lemma:Fourierseriesnew} that for all $x \in \big[-\frac{1}{2},\,\frac{1}{2}\big]$ it holds,
\begin{eqnarray*}
& &\frac{1}{N_1}\,\sum_{\ell\in \I_{N_1+2m_1}} {\mathrm e}^{- 2\pi {\mathrm i}\,w \ell/N_1}\,{\mathrm e}^{-2 \pi {\mathrm i}\, w\,x}\,\varphi_1\bigg(\!x +\frac{\ell}{N_1}\!\bigg) - {\hat \varphi}_1(w)\\
& & =\,\sum_{r \in {\mathbb Z}\setminus \{0\}} {\hat \varphi}_1(w + r N_1)\,{\mathrm e}^{2 \pi {\mathrm i}\,r N_1\,x}\,.
\end{eqnarray*}
Since by assumption ${\hat \varphi}_1(w) > 0$ for all $w \in \big[-\frac{N}{2a},\,\frac{N}{2a}\big]\subset \big[-\frac{N}{2},\,\frac{N}{2}\big]$, we have
\begin{eqnarray*}
& &\frac{1}{N_1 {\hat \varphi}_1(w)}\,\sum_{\ell\in \I_{N_1+2m_1}} {\mathrm e}^{- 2\pi {\mathrm i}\,w \ell/N_1}\,{\mathrm e}^{-2 \pi {\mathrm i}\, w\,x}\,\varphi_1\bigg(\!x +\frac{\ell}{N_1}\!\bigg) - 1\\
& & =\,\sum_{r \in {\mathbb Z}\setminus \{0\}} \frac{{\hat \varphi}_1(w + r N_1)}{{\hat \varphi}_1(w)}\,{\mathrm e}^{2 \pi {\mathrm i}\,r N_1\,x}\,.
\end{eqnarray*}
Multiplying the above equality by the exponential ${\mathrm e}^{2 \pi {\mathrm i}\, w\,x}$, this results in \eqref{eq:sumhatvarphi} and \eqref{eq:maxmax}. \qedsymbol
\medskip

We say that the window function $\varphi_1$ of the form \eqref{eq:firstwindow} is \emph{convenient for} NNFFT, if the \emph{general} $C(\mathbb T)$-\emph{error constant}
\begin{equation}
\label{eq:Esigma(varphi)}
E_{\sigma_1}(\varphi_1) \coloneqq \sup_{N\in\mathbb N} E_{\sigma_1,N}(\varphi_1)
\end{equation}
with
\begin{equation}
\label{eq:EsigmaN}
E_{\sigma_1,N}(\varphi_1) \coloneqq \max_{v \in [- N/2,\,  N/2]}  \bigg\| \sum_{r \in \mathbb Z \setminus \{0\}} \frac{{\hat \varphi}_1(v + r N_1)}{{\hat \varphi}_1(v)}\,{\mathrm e}^{2 \pi {\mathrm i}\,r N_1 \cdot}\bigg\|_{C(\mathbb T)}
\end{equation}
fulfills the condition
$E_{\sigma_1}(\varphi_1) \ll 1$
for conveniently chosen truncation parameter $m_1 \ge 2$ and oversampling factor $\sigma_1 > 1$. Obviously, the $C(\mathbb T)$-error constant \eqref{eq:esigma(varphi)} is a ``discrete'' version of the general $C(\mathbb T)$-error constant \eqref{eq:Esigma(varphi)} with the property
\begin{equation}
\label{eq:eleE}
e_{\sigma_1}(\varphi_1) \le E_{\sigma_1}(\varphi_1)\,.
\end{equation}
Thus, Corollary~\ref{Cor:Cor2.3} means that all complex exponentials ${\mathrm e}^{2 \pi {\mathrm i}\,w x}$ with $w\in \big[-\frac{N}{2a},\,\frac{N}{2a}\big]$ and $x\in \big[-\frac{1}{2},\,\frac{1}{2}\big]$ can be
uniformly approximated by short linear combinations of shifted window functions, cf. \cite[Theorem~2.10]{DuRo93}, if $\varphi_1$ is convenient for NNFFT.

\begin{Theorem}
\label{Thm:Esigma1(varphi1)}
Let $\sigma_1 >1$, $m_1 \in {\mathbb N} \setminus \{1\}$, and $N_1 = \sigma_1 N \in 2\mathbb N$ with $2 m_1 \ll N_1$ be given. Let $\varphi_1$ be the scaled version \eqref{eq:firstwindow} of $\omega_1 \in \Omega$.
Assume that the Fourier transform ${\hat \omega}_1$ fulfills the decay condition
$$
|{\hat \omega}_1(v)| \le \left\{ \begin{array}{ll} c_1 & \quad |v| \in \big[m_1 \big(1 - \frac{1}{2\sigma_1}\big),\,m_1 \big(1 + \frac{1}{2\sigma_1}\big)\big]\,,\\ [1ex]
 c_2 \,|v|^{-\mu} & \quad |v| \ge m_1 \big(1 + \frac{1}{2\sigma_1}\big)\,,
\end{array} \right.
$$
with certain constants $c_1 > 0$, $c_2>0$, and $\mu > 1$.

Then the general $C(\mathbb T)$-error constant $E_{\sigma_1}(\varphi_1)$ of the window function \eqref{eq:firstwindow} has the upper bound
\begin{equation}
\label{eq:boundEsigma(varphi)}
E_{\sigma_1}(\varphi_1) \le \frac{1}{{\hat \omega}_1\big(\frac{m_1}{2 \sigma_1}\big)} \Bigg[2 c_1 + \frac{2 c_2}{(\mu -1)\,m_1^{\mu}}\bigg(\!1 - \frac{1}{2\sigma_1}\!\bigg)^{1-\mu}\Bigg]\,.
\end{equation}
\end{Theorem}

\emph{Proof.}
By the scaling property of the Fourier transform, we have
$$
{\hat \varphi}_1(v) = \int_{\mathbb R} \varphi_1(t)\,{\mathrm e}^{-2 \pi {\mathrm i}v t}\,{\mathrm d}t = \frac{m}{N_1}\,{\hat \omega}_1\Big(\frac{m_1 v}{N_1}\Big)\,, \quad v \in \mathbb R\,.
$$
For all $v \in \big[-\frac{N}{2},\, \frac{N}{2}\big]$ and $r\in {\mathbb Z}\setminus \{0,\,\pm 1\}$, we obtain
$$
\bigg| \frac{m_1 v}{N_1} + m_1 r \bigg| \ge m_1 \bigg(\!2 - \frac{1}{2 \sigma_1}\!\bigg) > m_1 \bigg(\!1 + \frac{1}{2 \sigma_1}\!\bigg)
$$
and hence
$$
|{\hat \varphi_1}(v + r N_1)| = \frac{m_1}{N_1}\,\Big|{\hat \omega}_1\Big(\frac{m_1 v}{N_1} + m_1 r\Big)\Big| \le \frac{m_1\,c_2}{m_1^{\mu} N_1}\,\Big|\frac{v}{N_1} + r\Big|^{-\mu}\,.
$$
From \cite[Lemma~3.1]{PoTa20} it follows that for fixed $u =\frac{v}{N_1}\in \big[- \frac{1}{2\sigma_1},\, \frac{1}{2\sigma_1}\big]$,
$$
\sum_{r \in {\mathbb Z} \setminus \{0, \pm 1\}} |u + r|^{-\mu} \le \frac{2}{\mu - 1}\, \bigg(\!1 - \frac{1}{2\sigma_1}\!\bigg)^{1-\mu}\,.
$$
For all $v \in \big[- \frac{N}{2},\, \frac{N}{2}\big]$, we sustain
$$
|{\hat \varphi}_1(v \pm N_1)| = \frac{m_1}{N_1}\,\Big|{\hat \omega}_1 \Big(\frac{m_1 v}{N_1} \pm m_1 \Big)\Big| \le \frac{m_1}{N_1}\,c_1\,,
$$
since it holds
$$
\Big|\frac{m_1 v}{N_1} \pm m_1 \Big| \in \bigg[m_1 \bigg(\!1 - \frac{1}{2\sigma_1}\!\bigg),\,m_1 \bigg(\!1 + \frac{1}{2\sigma_1}\!\bigg)\bigg]\,.
$$
Thus, for each $v \in \big[-\frac{N}{2},\,\frac{N}{2}\big]$, we estimate the sum
\begin{eqnarray*}
\sum_{r \in {\mathbb Z} \setminus \{0\}} |{\hat \varphi}_1(v + r N_1)| &\le& \frac{m_1}{N_1}\,\bigg[\Big|{\hat \omega}_1\Big(\frac{m_1 v}{N_1} - m_1 \Big)\Big| + \Big|{\hat \omega}_1\Big(\frac{m_1 v}{N_1} + m_1 \Big)\Big|\\
& & + \,\sum_{k \in {\mathbb Z} \setminus \{0, \pm 1\}} \Big|{\hat \omega}_1\Big(\frac{m_1 v}{N_1} + m_1 r\Big)\Big|\bigg]\\
&\le& \frac{m_1}{N_1}\, \bigg[2 c_1 + \frac{c_2}{m_1^{\mu}}\,\sum_{r \in {\mathbb Z} \setminus \{0,\pm 1\}} \Big|\frac{v}{N_1} + r \Big|^{- \mu}\bigg]\\
&\le& \frac{m_1}{N_1}\,\bigg[2 c_1 + \frac{2 c_2}{(\mu - 1) \,m_1^{\mu}}\, \Big(1 -\frac{1}{2\sigma_1}\Big)^{1 -\mu}\bigg]
\end{eqnarray*}
such that
$$
\max_{v \in [-N/2,N/2]}\sum_{r \in {\mathbb Z} \setminus \{0\}} |{\hat \varphi}_1(v + r N_1)| \le \frac{m_1}{N_1}\,\Bigg[2 c_1 + \frac{2 c_2}{(\mu - 1) \,m_1^{\mu}}\, \bigg(\!1 -\frac{1}{2\sigma_1}\!\bigg)^{1 -\mu}\Bigg]\,.
$$
Now we determine the minimum of all positive values
$$
{\hat \varphi}_1(v) = \frac{m_1}{N_1}\, {\hat \omega}_1\Big(\frac{m_1 v}{N_1}\Big)\,, \quad v \in \bigg[\! - \frac{N}{2},\, \frac{N}{2}\!\bigg]\,.
$$
Since $\frac{m_1\,|v|}{N_1} \le \frac{m_1}{2 \sigma_1}$ for all $v \in \big[ - \frac{N}{2},\, \frac{N}{2}\big]$, we obtain
$$
\min_{v \in [-N/2,N/2]} {\hat \varphi}_1(v) = \frac{m_1}{N_1}\,\min_{v \in [-N/2,N/2]} {\hat \omega}_1\Big(\frac{m_1 v}{N_1}\Big)
= \frac{m_1}{N_1}\,{\hat \omega}_1\Big(\frac{m_1}{2 \sigma_1}\Big) = {\hat \varphi}_1\Big(\frac{N}{2}\Big) > 0\,.
$$
Thus, we see that the constant $E_{\sigma_1,N}(\varphi_1)$ can be estimated by an upper bound which depends
on $m_1$ and $\sigma_1$, but does not depend on $N$. We obtain
\begin{eqnarray*}
E_{\sigma_1,N}(\varphi_1) &\le& \frac{1}{{\hat \varphi}_1(N/2)}\,\max_{v \in [-N/2, N/2]} \sum_{r \in {\mathbb Z} \setminus \{0\}} |{\hat \varphi}_1(n + r N_1)|\\
&\le& \frac{1}{{\hat \omega}_1\big(\frac{m_1}{2\sigma_1}\big)}\,\Bigg[2 c_1 +
\frac{2 c_2}{(\mu - 1)\,m_1^{\mu}}\, \bigg(1 -\frac{1}{2\sigma_1}\bigg)^{1-\mu}\Bigg]\,.
\end{eqnarray*}
Consequently, the general $C(\mathbb T)$-error constant $E_{\sigma_1}(\varphi_1)$ has the upper bound \eqref{eq:boundEsigma(varphi)}. By \eqref{eq:eleE}, the expression \eqref{eq:boundEsigma(varphi)}
is also an upper bound of $C(\mathbb T)$-error constant $e_{\sigma_1}(\varphi_1)$. \qedsymbol
\medskip

Thus, \new{by means of these technical results} we obtain the following error estimate for the NNFFT.

\begin{Theorem}
\label{Thm:errorestimate}
Let the nonharmonic bandwidth $N\in \mathbb N$ with $N\gg 1$ be given.
Assume that $N_1 = \sigma_1 N \in 2 \mathbb N$ with $\sigma_1 > 1$. For fixed $m_1 \in \mathbb N \setminus \{1\}$ with $2 m_1 \ll N_1$, let  $N_2 = \sigma_2\,(N_1 + 2m_1)$ with $\sigma_2 > 1$. For $m_2 \in \mathbb N \setminus \{1\}$ with $2m_2 \le \big(1 - \frac{1}{\sigma_1}\big)\, N_2$, let
$\varphi_1$ and $\varphi_2$ be the window functions of the form \eqref{eq:twowindows}. Let $x_j \in \big[ - \frac{1}{2},\,\frac{1}{2}\big]$, $j \in \I_{M_2}$, be arbitrary spatial nodes
and let $f_k \in \mathbb C$, $k \in \I_{M_1}$, be arbitrary coefficients.
Further, let $a > 1$ be the constant \eqref{eq:consta}.

Then for a given exponential sum \eqref{eq:f(x)}  with arbitrary frequencies $v_k\in \big[- \frac{1}{2a},\,\frac{1}{2a}\big]$, $k \in \I_{M_1}$, the error of the $\mathrm{NNFFT}$ can be estimated by
\begin{align}
\label{eq:error_bound_nnfft}
	\max_{j\in \I_{M_2}} \bigg| f(x_j) - \frac{s_1(Nx_j)}{{\hat\varphi}_1(Nx_j)}\bigg| &\le \max_{x \in [-1/2,\,1/2]} \bigg| f(x) - \frac{s_1(Nx)}{{\hat\varphi}_1(Nx)}\bigg| \notag \\
	&\le \bigg[E_{\sigma_1}(\varphi_1) + \frac{a}{ {\hat \varphi}_1\big(\frac{N}{2}\big)}\,E_{\sigma_2}(\varphi_2)\bigg]\, \sum_{k\in \I_{M_1}} |f_k|\,,
\end{align}
where $E_{\sigma_j}(\varphi_j)$ for $j = 1,\,2,$ are the general $C(\mathbb T)$-error constants of the form \eqref{eq:Esigma(varphi)}.
\end{Theorem}

\new{\emph{Proof.}}
Now for arbitrary spatial nodes $x_j \in \big[-\frac{1}{2},\, \frac{1}{2}\big]$, $j\in \I_{M_2}$, we estimate the error of the NNFFT  in the form
$$
\max_{j\in \I_{M_2}} \bigg| f(x_j) - \frac{s_1(N \,x_j)}{{\hat \varphi}_1(Nx_j)}\bigg|\le \max_{j\in \I_{M_2}} \bigg| f(x_j) - \frac{s(Nx_j)}{{\hat \varphi}_1(N \,x_j)}\bigg| + \max_{j\in \I_{M_2}} \frac{|s(Nx_j) - s_1(Nx_j)|}{{\hat \varphi}_1(N \,x_j)}\,.
$$
At first we consider
$$
\max_{j\in \I_{M_2}} \bigg| f(x_j) - \frac{s(N x_j)}{{\hat \varphi}_1(Nx_j)}\bigg| \le \max_{x\in [-1/2,\,1/2]} \bigg| f(x) - \frac{s(Nx)}{{\hat \varphi}_1(Nx)}\bigg|\,.
$$
From \eqref{eq:hath(v)} and \eqref{eq:s(v):=} it follows that for all $x \in \mathbb R$,
\begin{eqnarray*}
	& &f(x) - \frac{s(Nx)}{{\hat \varphi}_1(N x)} = \frac{\hat h(Nx) - s(Nx)}{\hat \varphi_1(Nx)}\\
	& &=\, \sum_{k\in \I_{M_1}} f_k\,\Bigg[ {\mathrm e}^{-2\pi{\mathrm i}N v_k x} - \frac{1}{N_1\,{\hat \varphi}_1(Nx)} \sum_{\ell \in \I_{N_1+2m_1}}\varphi_1\bigg(\!\frac{\ell}{N_1}-v_k\!\bigg)\,{\mathrm e}^{-2\pi {\mathrm i}\, \ell x/\sigma_1}\Bigg]\,.
\end{eqnarray*}

Thus, by \eqref{eq:condvk}, \eqref{eq:maxmax}, and \eqref{eq:EsigmaN}, we obtain the estimate
\begin{align}
\label{eq:E1}
	\max_{x\in [-1/2,\,1/2]} \bigg| f(x) - \frac{s(Nx)}{{\hat \varphi}_1(Nx)}\bigg| \le E_{\sigma_1,N}(\varphi_1)\,\sum_{k\in \I_{M_1}} |f_k| \le E_{\sigma_1}(\varphi_1)\,\sum_{k\in \I_{M_1}} |f_k|\,.
\end{align}
Now we show that for $\varphi_2(t) \coloneqq \omega_2\big(\frac{N_2 t}{m_2}\big)$ and $N_2 = \sigma_2\,(N_1 + 2 m_1)$ it holds
\begin{equation}
	\label{eq:max|s(Nx)-s1(Nx)|}
	\max_{x\in [-1/2,\,1/2]} |s(N \,x) - s_1(Nx)| \le  E_{\sigma_2}(\varphi_2)\,\sum_{\ell\in \I_{N_1+ 2m_1}} |g_{\ell}|\,.
\end{equation}
By construction, the functions $s$ and $s_1$ can be represented in the form
\begin{eqnarray*}
	s(N x) &=& \sum_{\ell \in \I_{N_1 + 2m_1}} g_{\ell}\,{\mathrm e}^{-2 \pi {\mathrm i}\ell x/\sigma_1}\,,\\
	s_1(N x) &=& \sum_{s \in \I_{N_2}} h_s\, {\tilde \varphi}_2^{(1)}\bigg(\!\frac{x}{\sigma_1} - \frac{s}{N_2}\!\bigg)\,, \quad x \in \mathbb R\,,
\end{eqnarray*}
where ${\tilde \varphi}_2^{(1)}$ denotes the 1-periodization of the second window function $\varphi_2$ and
$$
h_s \coloneqq \frac{1}{N_2} \sum_{\ell \in \I_{N_1 + 2m_1}} \frac{g_{\ell}}{{\hat \varphi}_2(\ell)}\,{\mathrm e}^{- 2 \pi {\mathrm i} \ell s /N_2}\,.
$$
Substituting $t = \frac{x}{\sigma_1}$, it follows that
\begin{eqnarray*}
	s(N_1 t) &=& \sum_{\ell \in \I_{N_1 + 2m_1}} g_{\ell}\,{\mathrm e}^{-2 \pi {\mathrm i}\ell t}\,,\\
	s_1(N_1 t) &=& \sum_{s \in \I_{N_2}} h_s\, {\tilde \varphi}_2^{(1)}\bigg(\!t - \frac{s}{N_2}\!\bigg)\,, \quad t \in \mathbb R\,,
\end{eqnarray*}
are 1-periodic functions. By \cite[Lemma~2.3]{PoTa20}, we conclude
$$
\max_{t\in [-1/2,\,1/2]} |s(N_1 t) - s_1(N_1 t)| \le\, e_{\sigma_2}(\varphi_2)\,\sum_{\ell\in \I_{N_1+ 2m_1}} |g_{\ell}| \le  E_{\sigma_2}(\varphi_2)\,\sum_{\ell\in \I_{N_1+ 2m_1}} |g_{\ell}|\,,
$$
where the general $C(\mathbb T)$-error constant $E_{\sigma_2}(\varphi_2)$ defined similar to \eqref{eq:Esigma(varphi)} has an analogous property \eqref{eq:eleE}.
Since $x = \sigma_1\,t$, we obtain that
\begin{align}
\label{eq:E2}
	\max_{t\in [-1/2,\,1/2]} |s(N_1 t) - s_1(N_1 t)| &= \max_{x \in [-\sigma_1/2, \sigma_1/2]} |s(N x) - s_1(N x)|
	\notag \\
	&\le E_{\sigma_2}(\varphi_2)\,\sum_{\ell\in \I_{N_1+ 2m_1}} |g_{\ell}|\,,
\end{align}
such that \eqref{eq:max|s(Nx)-s1(Nx)|} is shown. Note that for $x \in \big[- \frac{1}{2}, \,\frac{1}{2}\big]$ it holds
$$
s_1(N x) = \sum_{s \in \I_{N_2}''(x)} h_s\,\varphi_2\bigg(\!\frac{x}{\sigma_1} - \frac{s}{N_2}\!\bigg)
$$
with the index set
$$
\I_{N_2}''(x) \coloneqq \bigg\{ s \in \I_{N_2}:\, \bigg|\frac{s}{N_2} - \frac{x}{\sigma_1}\bigg| < \frac{m_2}{N_2}\bigg\}\,.
$$
Further, by \eqref{eq:consta} and \eqref{eq:gell} it holds
\begin{eqnarray*}
	\sum_{\ell\in \I_{N_1+ 2m_1}} |g_{\ell}| &\le& \frac{1}{N_1}\,\sum_{\ell \in \I_{N_1+2m_1}} \sum_{k\in \I_{M_1}} |f_k| \cdot 1 \le \frac{N_1 + 2m_1}{N_1} \,\sum_{k\in \I_{M_1}} |f_k|
	= a \,\sum_{k\in \I_{M_1}} |f_k|\,.
\end{eqnarray*}
\new{Combining this with \eqref{eq:E1} and \eqref{eq:E2} completes the proof.}
\qedsymbol
\medskip

\new{Now it merely remains to estimate the general $C(\mathbb T)$-error constants $E_{\sigma_j}(\varphi_j)$ for $j = 1,\,2,$ and ${\hat \varphi}_1\big(\frac{N}{2}\big)$ in \eqref{eq:error_bound_nnfft} for specific window functions.}

\section{Error of NNFFT with sinh-type window functions}
\label{Sec:NNFFTsinhwindow}

\new{In this section we specify the result in Theorem~\ref{Thm:errorestimate} for the NNFFT with two $\sinh$-type window functions.}

Let \mbox{$N\in \mathbb N$} with \mbox{$N\gg 1$} be the fixed nonharmonic bandwidth.
Let \mbox{$\sigma_1,\, \sigma_2 \in \big[\frac{5}{4},\,2\big]$} be given oversampling factors.
Further let \mbox{$N_1 = \sigma_1 N \in 2 \mathbb N$}, \mbox{$m_1 \in \mathbb N \setminus \{1\}$} with \mbox{$2 m_1 \ll N_1$}, and \mbox{$N_2 = \sigma_2\,(N_1 + 2m_1) = \sigma_1 \sigma_2 a\, N \in 2 \mathbb N$}
be given, where \mbox{$a>1$} denotes the constant \eqref{eq:consta}.
Let \mbox{$m_2 \in \mathbb N \setminus \{1\}$} with \mbox{$2m_2 \le \big(1 - \frac{1}{\sigma_1}\big)\, N_2$} be given as well.

For \mbox{$j = 1,\,2$}, we consider the functions
\begin{equation*}
\omega_{\sinh,j}(x) \coloneqq \left\{\begin{array}{ll} \frac{1}{\sinh \beta_j}\, \sinh\big(\beta_j\,\sqrt{1 - x^2}\big) & \quad x\in [-1,\,1]\,,\\
0 & \quad x \in {\mathbb R} \setminus [-1,\,1]
\end{array} \right.
\end{equation*}
with the shape parameter
$$
\beta_j \coloneqq 2 \pi m_j \bigg(\!1 - \frac{1}{2\sigma_j}\!\bigg)\,.
$$
As shown in Example~\ref{Example:omega}, both functions belong to the set $\Omega$. By scaling, for \mbox{$j = 1,\,2$}, we introduce the $\sinh$-\emph{type window functions}
\begin{equation}
\label{eq:varphisinhj}
\varphi_{\sinh,j}(t) \coloneqq \omega_{\sinh,j}\bigg(\!\frac{N_j t}{m_j}\!\bigg)\,, \quad t \in \mathbb R\,.
\end{equation}

Now we show that the error of the NNFFT with two $\sinh$-type window functions \eqref{eq:varphisinhj} has exponential decay with respect to the truncation parameters $m_1$ and $m_2$.

\begin{Theorem}
\label{Thm:errorNNFFTsinhwindow}
Let the nonharmonic bandwidth \mbox{$N\in\mathbb N$} with \mbox{$N\gg 1$} be given.
Further let \mbox{$N_1 = \sigma_1 N \in 2 \mathbb N$} with \mbox{$\sigma_1 \in \big[\frac{5}{4},\,2\big]$} be given.
For fixed \mbox{$m_1 \in \mathbb N \setminus \{1\}$} with \mbox{$2 m_1 \ll N_1$}, let
\mbox{$N_2 = \sigma_2\,(N_1 + 2m_1)\in 2 \mathbb N$} with
\mbox{$\sigma_2 \in \big[\frac{5}{4},\,2\big]$}.
For \mbox{$m_2 \in \mathbb N \setminus \{1\}$} with
\mbox{$2m_2 \le \big(1 - \frac{1}{\sigma_1}\big)\, N_2$}, let
$\varphi_{\sinh,1}$ and $\varphi_{\sinh,2}$ be the $\sinh$-type window functions \eqref{eq:varphisinhj}.
Assume that \mbox{$m_2 \ge m_1$}.
Let \mbox{$x_j \in \big[ - \frac{1}{2},\,\frac{1}{2}\big]$}, \mbox{$j \in \I_{M_2}$},
be arbitrary spatial nodes and let
\mbox{$f_k \in \mathbb C$}, \mbox{$k \in \I_{M_1}$}, be arbitrary coefficients.
Let \mbox{$a>1$} be the constant \eqref{eq:consta}.

Then for the exponential sum \eqref{eq:f(x)} with arbitrary frequencies
\mbox{$v_k\in \big[- \frac{1}{2a},\,\frac{1}{2a}\big]$}, \mbox{$k \in \I_{M_1}$},
the error of the $\mathrm{NNFFT}$ with the $\sinh$-type window functions \eqref{eq:varphisinhj} can be estimated in the form
$$
\max_{j\in \I_{M_2}} \bigg| f(x_j) - \frac{s_1(Nx_j)}{{\hat\varphi_{\sinh,1}(Nx_j)}}\bigg| \le \max_{x \in [-1/2,\,1/2]} \bigg| f(x) - \frac{s_1(Nx)}{{\hat\varphi_{\sinh,1}(Nx)}}\bigg| \le E(\varphi_{\sinh})\,\sum_{k\in \I_{M_1}} |f_k|
$$
with the constant
\begin{align}
\label{eq:error_bound}
E(\varphi_{\sinh}) &\coloneqq (24 m_1^{3/2} + 10)\,{\mathrm e}^{- 2 \pi m_1 \sqrt{1- 1/\sigma_1}} \notag \\
& +\, (24 m_2^{3/2} + 10)\,\frac{2N_1 + 4m_1}{\sqrt{2 m_1 \pi}}\,{\mathrm e}^{2 \pi m_1\,(1- \sqrt{1 - 1/\sigma_1}- 1/(2\sigma_1))}\,{\mathrm e}^{- 2 \pi m_2 \sqrt{1- 1/\sigma_2}}\,.
\end{align}
\end{Theorem}

\new{\emph{Proof.}}
\new{By Theorem~\ref{Thm:errorestimate} we have to estimate the general $C(\mathbb T)$-error constants \mbox{$E_{\sigma_j}(\varphi_j)$}, \mbox{$j = 1,\,2,$} and ${\hat \varphi}_1\big(\frac{N}{2}\big)$ in \eqref{eq:error_bound_nnfft} for the $\sinh$-type window functions \eqref{eq:varphisinhj}.}

Applying Theorem~\ref{Thm:Esigma1(varphi1)}, we obtain by the same technique as in \cite[Theorem~5.6]{PoTa20} that
\begin{equation}
	\label{eq:Esigmaj(varphisinhj)}
	E_{\sigma_j}\big(\varphi_{\sinh,j}\big) \le (24 m_j^{3/2} + 10)\,{\mathrm e}^{- 2 \pi m_j \sqrt{1- 1/\sigma_j}}\,, \quad j=1,\,2.
\end{equation}
Now we estimate ${\hat \varphi}_{\sinh,1}\big(\frac{N}{2}\big)$. Using the scaling property of the Fourier transform, by \eqref{eq:hatomegasinh} we obtain
\begin{eqnarray*}
	{\hat \varphi}_{\sinh,1}\bigg(\!\frac{N}{2}\!\bigg) &=& \frac{m_1}{N_1}\,{\hat \omega}_{\sinh,1}\bigg(\!\frac{m_1 N}{2 N_1}\!\bigg) = \frac{m_1}{N_1}\,{\hat \omega}_{\sinh,1}\bigg(\!\frac{m_1}{2 \sigma_1}\!\bigg)\\
	&=& \frac{\pi m_1 \beta_1}{N_1 \sinh \beta_1}\, \bigg(\!\beta_1^2 - \frac{\pi^2 m_1^2}{\sigma_1^2}\!\bigg)^{-1/2} \,I_1\Bigg(\!\sqrt{\beta_1^2 - \frac{\pi^2 m_1^2}{\sigma_1^2}}\,\Bigg)\\
	&=& \frac{m_1 \pi}{N_1 \sinh \beta_1}\,\bigg(\!1 - \frac{1}{2 \sigma_1}\!\bigg)\bigg(\!1 - \frac{1}{\sigma_1}\!\bigg)^{-1/2}\, I_1\bigg(\!2 \pi m_1 \sqrt{1-\frac{1}{\sigma_1}}\,\bigg)\,,
\end{eqnarray*}
where we have used the equality
$$
\bigg(\!\beta_1^2 - \frac{\pi^2 m_1^2}{\sigma_1^2}\!\bigg)^{1/2} = 2\pi m_1 \Bigg(\! \bigg(\!1- \frac{1}{2\sigma_1}\!\bigg)^2 - \frac{1}{4 \sigma_1^2}\Bigg)^{1/2} = 2 \pi m_1 \sqrt{1 - \frac{1}{\sigma_1}}\,.
$$
From $m_1 \ge 2$ and $\sigma_1 \ge \frac{5}{4}$, it follows that
$$
2 \pi m_1\, \sqrt{1 - \frac{1}{\sigma_1}} \ge 4 \pi \, \sqrt{1 - \frac{1}{\sigma_1}} \ge x_0 \coloneqq \frac{4 \pi}{\sqrt 5}\,.
$$
By the inequality for the modified Bessel function $I_1$ (see \cite[Lemma~3.3]{PoTa20}) it holds
$$
I_1(x) \ge \sqrt x_0\,{\mathrm e}^{-x_0}\,I_1(x_0)\, x^{-1/2}\,{\mathrm e}^x > \frac{2}{5}\, x^{-1/2}\,{\mathrm e}^x\,, \quad x \ge x_0\,.
$$
Thus, we obtain
$$
{\hat \varphi}_{\sinh,1}\bigg(\!\frac{N}{2}\!\bigg) \ge \frac{\sqrt{2 m_1 \pi}}{5 N_1\,\sinh \beta_1}\,\bigg(\!1- \frac{1}{2\sigma_1}\!\bigg)\bigg(\!1- \frac{1}{\sigma_1}\!\bigg)^{-3/4}\,{\mathrm e}^{2 \pi m_1 \sqrt{1 - 1/\sigma_1}}\,.
$$
By the simple inequality
$$
\sinh \beta_1 < \frac{1}{2}\,{\mathrm e}^{\beta_1} = \frac{1}{2}\,{\mathrm e}^{2 \pi m_1 (1- 1/(2\sigma_1))}\,,
$$
we conclude that
$$
{\hat \varphi}_{\sinh,1}\bigg(\!\frac{N}{2}\!\bigg) \ge \frac{2 \sqrt{2 m_1 \pi}}{5 N_1}\,\bigg(\!1- \frac{1}{2\sigma_1}\!\bigg)\bigg(\!1- \frac{1}{\sigma_1}\!\bigg)^{-3/4}\,{\mathrm e}^{2 \pi m_1 (\sqrt{1 - 1/\sigma_1} - 1 + 1/(2\sigma_1))}
$$
and hence
\begin{equation}
	\label{eq:a/hatvarphisinh1}
	\frac{a}{{\hat \varphi}_{\sinh,1}(N/2)} \le \frac{5 N_1\,a}{2 \sqrt{2 m_1 \pi}}\,\bigg(\!1- \frac{1}{2\sigma_1}\!\bigg)^{-1}\bigg(\!1- \frac{1}{\sigma_1}\!\bigg)^{3/4}\,{\mathrm e}^{-2 \pi m_1 (\sqrt{1 - 1/\sigma_1} - 1 + 1/(2\sigma_1))}\,.
\end{equation}
Applying Theorem~\ref{Thm:errorestimate}, we estimate the error of the NNFFT with two $\sinh$-type window functions \eqref{eq:varphisinhj}.
By \eqref{eq:Esigmaj(varphisinhj)} and \eqref{eq:a/hatvarphisinh1} we obtain  the inequality
\begin{eqnarray*}
	E_{\sigma_1}(\varphi_{\sinh,1}) &+& \frac{a}{{\hat \varphi}_{\sinh,1}(N/2)}\,E_{\sigma_2}(\varphi_{\sinh,2}) \le (24 m_1^{3/2} + 10)\,{\mathrm e}^{- 2 \pi m_1 \sqrt{1- 1/\sigma_1}}\\
	&+&  \,(24 m_2^{3/2} + 10)\, \frac{2 N_1 a}{\sqrt{2 m_1 \pi}}\,{\mathrm e}^{2 \pi m_1 (1 - \sqrt{1 - 1/\sigma_1} - 1/(2\sigma_1))}\,{\mathrm e}^{- 2 \pi m_2 \sqrt{1- 1/\sigma_2}}\,,
\end{eqnarray*}
since it holds
$$
\frac{5}{2}\,\bigg(\!1- \frac{1}{2\sigma_1}\!\bigg)^{-1}\bigg(\!1- \frac{1}{\sigma_1}\!\bigg)^{3/4} \le \frac{5}{3}\sqrt{2}<2
\,, \quad \sigma_1 \in \big[\tfrac{5}{4},\, 2 \big]\,.
$$
\new{This completes the proof.}
\qedsymbol
\medskip

\begin{Example}
Now we visualize the result of Theorem~\ref{Thm:errorNNFFTsinhwindow}.
\new{To this end, we fix \mbox{$N=1200$} and consider \mbox{$m_1\in \{2, \ldots, 8\}$} and \mbox{$\sigma_1\in \{1.25,1.5,2\}$}.
In Figure~\ref{fig:error_bound} the error bound~\eqref{eq:error_bound} is depicted for several choices of \mbox{$m_2\geq m_1$} and \mbox{$\sigma_2\geq \sigma_1$}.
Clearly, the error bounds~\eqref{eq:error_bound} decrease for increasing truncation parameters and oversampling factors, respectively.
Moreover, we recognize that the results get better when choosing \mbox{$\sigma_2 > \sigma_1$}, cf.~Figure \ref{fig:error_bound}(c), and are best for \mbox{$m_2 > m_1$}, cf.~Figure \ref{fig:error_bound} (a).
Besides, we remark that choices \mbox{$m_2 < m_1$} or \mbox{$\sigma_2 < \sigma_1$} produce the same results as in the equality setting such that we omitted these tests.}

\new{Therefore, we recommend the use of truncation parameters \mbox{$m_2 > m_1$} and oversampling factors \mbox{$\sigma_2\geq \sigma_1$}.
For the choice of $m_1$ and $\sigma_1$, we refer to previous works concerning the NFFT, e.\,g.~\cite{PoTa20,PoTa20a}.}

\new{Additionally, we aim to compare these theoretical bounds with the errors obtained by the NNFFT.}
For this purpose, we introduce the relative error
\begin{equation}
\label{eq:relative_nnfft_error}
\bigg(\sum_{k \in \I_{M_1}} |f_k|\bigg)^{-1} \max_{j \in \I_{M_2}} \bigg| f(x_j) - \frac{s_1(N x_j)}{{\hat \varphi}_{\sinh,1}(N x_j)}\bigg|\,,
\end{equation}
since by Theorem~\ref{Thm:errorNNFFTsinhwindow} it holds
$$
\bigg(\sum_{k \in \I_{M_1}} |f_k|\bigg)^{-1} \max_{j \in \I_{M_2}} \bigg| f(x_j) - \frac{s_1(N x_j)}{{\hat \varphi}_{\sinh,1}(N x_j)}\bigg| \le E(\varphi_{\sinh})\,.
$$
Thus, we choose random nodes \mbox{$x_j \in \big[-\frac 12, \frac 12\big]$}, \mbox{$j \in \I_{M_2}$}, and \mbox{$v_k \in [- \frac{1}{2a}, \frac{1}{2a}]$}, \mbox{$k \in \I_{M_1}$}, with \mbox{$a = 1 + \frac{2m_1}{N_1}$},
as well as random coefficients \mbox{$f_k\in\mathbb C$}, \mbox{$k \in \I_{M_1}$}, and compute the values
\eqref{eq:f(xj)} once directly and once rapidly using the NNFFT.
Due to the randomness of the given data, this test is repeated one hundred times and afterwards the maximum error over all repetitions is computed.
The errors \eqref{eq:relative_nnfft_error} for the parameter choice \mbox{$M_1=2400$} and \mbox{$M_2=1600$} are displayed in Figure~\ref{fig:error_bound}~(b).

\new{Unfortunately, the current release NFFT 3.5.3 of the software package \cite{nfft3} is not yet designed for the use of parameters \mbox{$m_1\neq m_2$} and \mbox{$\sigma_1\neq \sigma_2$}.
Therefore, we can only handle the setting $m_1 = m_2$ and $\sigma_1 = \sigma_2$ in Figure~\ref{fig:error_bound}~(b).
Moreover, the $\sinh$-type window function is currently not implemented in the software package \cite{nfft3}.
Thus, we use two standard window functions, namely the Kaiser--Bessel window functions, since it was shown in \cite{PoTa20} that those are very much related.
Since the results in Figure~\ref{fig:error_bound} show great promise, these features might be part of future releases.}
\begin{figure}[ht]
	\centering
	\captionsetup[subfigure]{justification=centering}
	\begin{subfigure}[t]{0.32\textwidth}
		\includegraphics[width=\textwidth]{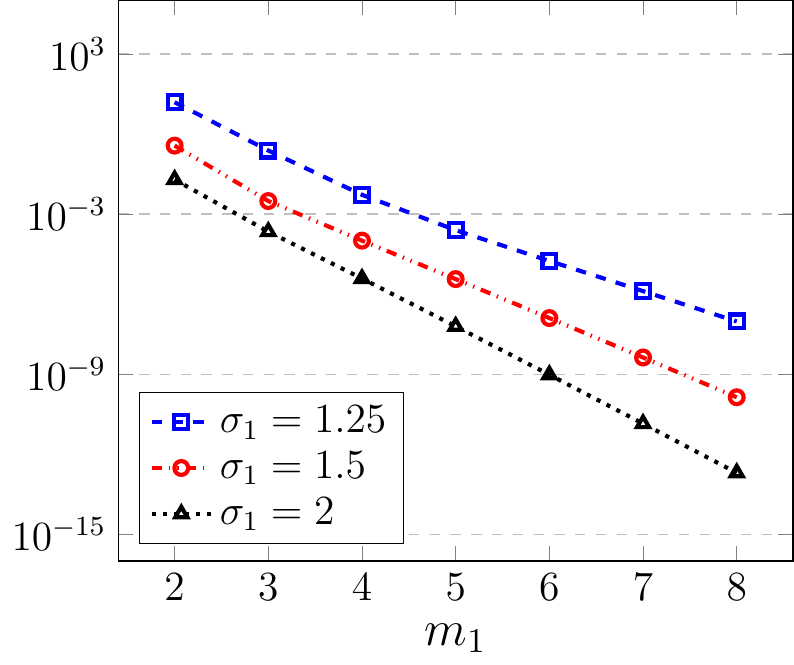}
		\caption{\new{$m_2=2m_1$ and $\sigma_2=\sigma_1$}}
	\end{subfigure}
	\begin{subfigure}[t]{0.32\textwidth}
		\includegraphics[width=\textwidth]{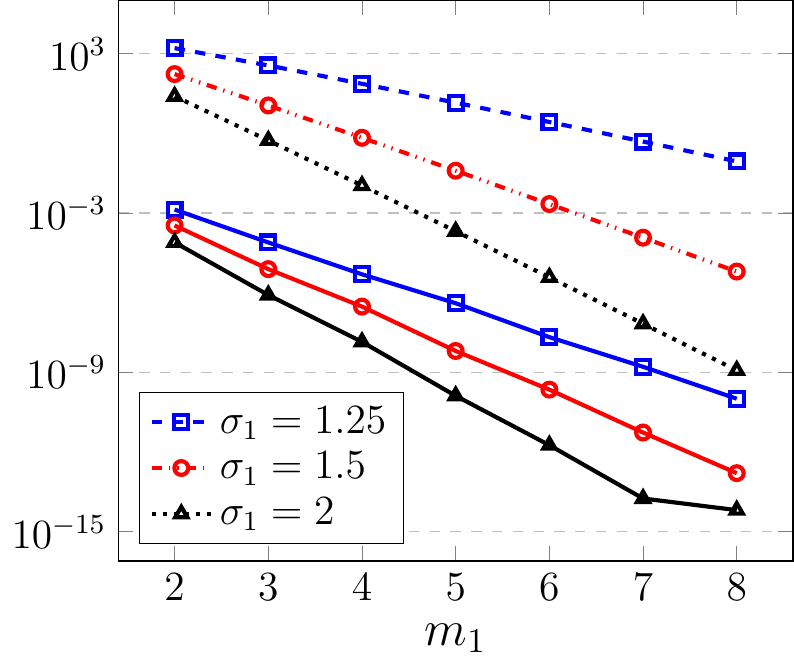}
		\caption{\new{$m_2=m_1$ and $\sigma_2=\sigma_1$}}
	\end{subfigure}
	\begin{subfigure}[t]{0.32\textwidth}
		\includegraphics[width=\textwidth]{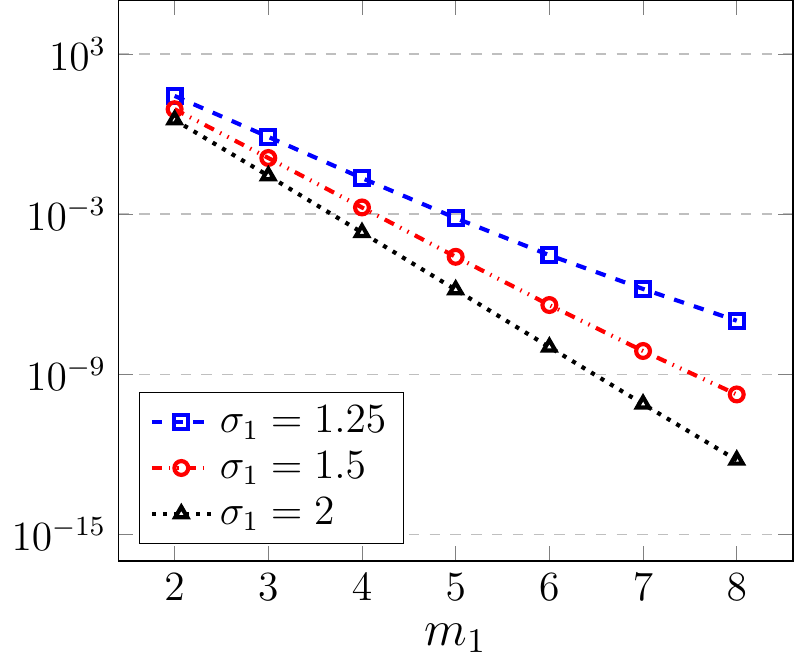}
		\caption{\new{$m_2=m_1$ and $\sigma_2=2\sigma_1$}}
	\end{subfigure}
	\caption{Error bound \eqref{eq:error_bound} (dashed) for the NNFFT with sinh-type window functions for \mbox{$N=1200$}, \mbox{$m_1\in\{2,\dots,8\}$} and $\sigma_1\in \{1.25,1.5,2\}$.
	\new{Part (b) additionally depicts the relative error \eqref{eq:relative_nnfft_error} (solid) using Kaiser-Bessel window functions.}
		\label{fig:error_bound}}
\end{figure}
\end{Example}

\section{Approximation of sinc function by exponential sums}
\label{Sec:Approxsinc}

\new{Since we aim to present an interesting signal processing application of the NNFFT in the last Section~\ref{Sec:discreteSinc}, we now study the approximation of the function $\mathrm{sinc}(N \pi x)$, \mbox{$x \in [-1,\,1]$}, by an exponential sum \eqref{eq:f(x)}.}

In \cite{BM02} the exponential sum \eqref{eq:f(x)} is used for a local approximation of a bandlimited function $F$ of the form
\begin{equation}
\label{eq:F(x)}
F(x) \coloneqq \int_{-1/2}^{1/2} w(t)\, {\mathrm e}^{-2 \pi {\mathrm i}N t x}\,{\mathrm d}t\,, \quad x \in \mathbb R\,,
\end{equation}
where $w:\,[-\frac{1}{2},\,\frac{1}{2}]\to [0,\,\infty)$ is an integrable function with \mbox{$\int_{-1/2}^{1/2} w(t)\,{\mathrm d}t >0$}.
By the substitution
\begin{equation*}
	F(x) = \frac{1}{N} \int_{-N/2}^{N/2} w\Big(\!-\frac{s}{N}\Big)\,{\mathrm e}^{2 \pi {\mathrm i}\, s x}\,{\mathrm d}s\,,
\end{equation*}
we recognize that the Fourier transform of \eqref{eq:F(x)} is supported on \mbox{$\big[-\frac{N}{2},\,\frac{N}{2}\big]$}, i.\,e., the function \eqref{eq:F(x)} is bandlimited with bandwidth $N$.
For instance, for \mbox{$w(t) \coloneqq 1$}, \mbox{$t \in \big[-\frac{1}{2},\,\frac{1}{2}\big]$}, we obtain the famous bandlimited $\mathrm{sinc}$ function
\begin{equation}
\label{eq:sinc}
F(x) = \mathrm{sinc}(\pi N x) \coloneqq \left\{ \begin{array}{ll} \frac{\sin (\pi N x)}{\pi N x} & \quad x \in \mathbb R \setminus \{0\}\,,\\
1 & \quad x = 0\,.
\end{array} \right.
\end{equation}
Now we show that the bandlimited $\mathrm{sinc}$ function \eqref{eq:sinc} can be uniformly approximated on the interval \mbox{$[-1,\,1]$} by an exponential sum \eqref{eq:f(x)}. We start with the uniform approximation
of the $\mathrm{sinc}$ function on the interval \mbox{$\big[-\frac{1}{2},\,\frac{1}{2}\big]$}.

\begin{Theorem}
	Let \mbox{$\varepsilon > 0$} be a given target
	accuracy. 
	
	Then for sufficiently large \mbox{$n\in \mathbb N$} with \mbox{$n \ge 2N$}, there exist constants \mbox{$w_j > 0$} and frequencies \mbox{$v_j\in \big(- \frac{1}{2},\,\frac{1}{2}\big)$}, \mbox{$j = 1,\ldots,n$}, such that for
	all \mbox{$x\in \big[-\frac{1}{2},\,\frac{1}{2}\big]$},
	\begin{equation}
	\label{eq:approx_nndft}
	\bigg| \mathrm{sinc}(\pi N x) - \sum_{j=1}^n w_j\, {\mathrm e}^{-2\pi {\mathrm i}Nv_j x}\bigg| \le \varepsilon\,.
	\end{equation}
\end{Theorem}

\emph{Proof}. This result is a simple consequence of \cite[Theorem~6.1]{BM02}.
Introducing \mbox{$\nu \coloneqq \frac{N}{n} \le \frac{1}{2}$},
we obtain by substitution \mbox{$\tau \coloneqq - \frac{t}{2 \nu}$} that
$$
\mathrm{sinc}(\pi N x) =  \int_{-1/2}^{1/2} {\mathrm e}^{-2 \pi {\mathrm i}N \tau x}\,{\mathrm d}\tau = \frac{1}{2\nu}\,\int_{-\nu}^{\nu} {\mathrm e}^{{\mathrm i} \pi n t x}\,{\mathrm d}t\,.
$$
Setting $y\coloneqq n x \in \big[-\frac{n}{2},\,\frac{n}{2}\big]$, we have
$$
\mathrm{sinc}(\pi \nu y) = \frac{1}{2 \nu}\,\int_{-\nu}^{\nu} {\mathrm e}^{{\mathrm i} \pi  t y}\,{\mathrm d}t\,.
$$
Then from  \cite[Theorem~6.1]{BM02} (with $d = \frac{1}{2}$), it follows the existence of \mbox{$w_j > 0$} and \mbox{$\Theta_j \in (-\nu,\, \nu)$}, \mbox{$j = 1, \ldots, n$}, such that for all \mbox{$y \in \big[- \frac{n}{2}-1,\, \frac{n}{2} + 1\big]$},
$$
\bigg|\frac{1}{2\nu}\, \int_{-\nu}^{\nu} \sigma(t)\,{\mathrm e}^{{\mathrm i} \pi  t y}\,{\mathrm d}t - \sum_{j=1}^n w_j\, {\mathrm e}^{\pi {\mathrm i}\,\Theta_j y}\bigg| \le \varepsilon\,.
$$
Hence, for all $x = \frac{y}{n}\in \big[-\frac{1}{2},\,\frac{1}{2}\big]$, we conclude that for \mbox{$v_j \coloneqq - \frac{\Theta_j}{2 \nu} \in \big(- \frac{1}{2},\, \frac{1}{2}\big)$}, \mbox{$j=1, \ldots,n$},
\begin{equation*}
\bigg| \frac{1}{2\sigma} \int_{-\nu}^{\nu} {\mathrm e}^{{\mathrm i} \pi n t x}\,{\mathrm d}t - \sum_{j=1}^n w_j\, {\mathrm e}^{\pi {\mathrm i}\,n \Theta_j x}\bigg|
= \, \bigg|\mathrm{sinc}(\pi N x)  -  \sum_{j=1}^n w_j\, {\mathrm e}^{-2\pi {\mathrm i}Nv_j x}\bigg| \le \varepsilon\,.
\end{equation*}
This completes the proof. \qedsymbol
\medskip

Substituting the variable \mbox{$x = \frac{t}{2}$}, \mbox{$t \in [-1,\,1]$},
the frequencies \mbox{$v_j=\frac{z_j}{2}$}, \mbox{$z_j\in(-1,\,1)$}, and
replacing the bandwidth $N$ in \eqref{eq:approx_nndft} by $2N$, we obtain the following uniform
approximation of the $\mathrm{sinc}$ function \eqref{eq:sinc} on the interval \mbox{$[-1,\,1]$} (after denoting $t$ by $x$ and $z_j$ by $v_j$ again):

\begin{corollary}
Let \mbox{$\varepsilon > 0$} be a given target accuracy.
	
Then for sufficiently large \mbox{$n\in \mathbb N$} with \mbox{$n \ge 4N$}, there exist constants \mbox{$w_j > 0$} and frequencies \mbox{$v_j\in (-1,\,1)$}, \mbox{$j = 1,\ldots,n$}, such that \eqref{eq:approx_nndft} holds for all \mbox{$x\in [-1,\,1]$}, i.\,e.,
\begin{equation*}
	\bigg| \mathrm{sinc}(\pi N x) - \sum_{j=1}^n w_j\, {\mathrm e}^{-\pi {\mathrm i}Nv_j x}\bigg| \le \varepsilon \,, \quad x\in [-1,\,1]\,.
\end{equation*}
\end{corollary}

\new{In practice, we simplify the approximation procedure of the function $\mathrm{sinc}(N \pi x)$.
Since for fixed \mbox{$N \in \mathbb N$}, it holds
$$
\mathrm{sinc} (N \pi x) = \frac{1}{2}\, \int_{-1}^1 {\mathrm e}^{- \pi {\mathrm i} N t x}\,{\mathrm d}t\,, \quad x \in \mathbb R\,,
$$
the approximation on the interval \mbox{$[-1, 1]$} }can efficiently be realized by means of the Clenshaw--Curtis quadrature (see \cite[pp.~143--153]{Tref13} or
\cite[pp.~357--364]{PlPoStTa18}).
Using this procedure for the integrand \mbox{$\frac{1}{2}\,{\mathrm e}^{- \pi {\mathrm i} N t x}$}, \mbox{$t \in [-1,\,1]$}, with fixed parameter \mbox{$x \in [-1,\,1]$}, the Chebyshev points
\mbox{$z_k =  \cos \frac{k \pi}{n}\in [-1, 1]$}, \mbox{$k = 0,\ldots, n$}, and the positive coefficients
\begin{equation}
\label{eq:wk=}
w_k = \left\{ \begin{array}{ll}
\frac{1}{n}\,\varepsilon_n(k)^2 \sum_{j=0}^{n/2} \varepsilon_n(2 j)^2 \frac{2}{1 - 4 j^2}\,\cos \frac{2jk \pi}{n} & \quad n \in 2 \mathbb N\,,\\
\frac{1}{n}\,\varepsilon_n(k)^2 \sum_{j=0}^{(n-1)/2} \varepsilon_n(2 j)^2 \frac{2}{1 - 4 j^2}\,\cos \frac{2jk \pi}{n} & \quad n \in 2 \mathbb N + 1\,,\\
\end{array} \right.
\end{equation}
with $\varepsilon_n(0) = \varepsilon_n(n) \coloneqq \frac{\sqrt 2}{2}$ and $\varepsilon_n(j) \coloneqq 1$, $j = 1, \ldots, n-1$ (see \cite[p.~359]{PlPoStTa18}), we obtain
$$
\mathrm{sinc}(N \pi x) = \frac{1}{2}\, \int_{-1}^1 {\mathrm e}^{- \pi {\mathrm i} N t x}\, {\mathrm d}t \approx \sum_{k = 0}^n w_k \, {\mathrm e}^{- \pi {\mathrm i} N z_k x}\,.
$$
Further the coefficients fulfill the condition (see \cite[p.~359]{PlPoStTa18})
\begin{equation}
\label{eq:sumwk=1}
\sum_{k=0}^n w_k = 1\,.
\end{equation}
Then we receive the following error estimate.

\begin{Theorem}
\label{Thm:errorapproxsinc}
Let $N\in\mathbb N$, $n = \nu N$ be given.
Let \mbox{$z_k =  \cos \frac{k \pi}{n}\in [-1, 1]$}, be the Chebyshev points,
let $w_k$, \mbox{$k = 0,\ldots, n$}, denote the coefficients \eqref{eq:wk=}, and set
\new{\mbox{$C \coloneqq \frac{\pi\,({\mathrm e}^2-1)}{2\,{\mathrm e}}$}}.

Then for all $x \in [-1, \,1]$, the approximation error of $\mathrm{sinc}(N \pi x)$ can be estimated in the form
\begin{equation}
\label{eq:estimateClenshawCurtis2}
\bigg|\mathrm{sinc}(N \pi x) - \sum_{k=0}^n w_k \,{\mathrm e}^{- \pi {\mathrm i} N z_k x}\bigg| \le \new{\frac{36\,(1+\mathrm e^{-2CN})}{35\,({\mathrm e}^2 -1)} \,\mathrm e^{-N(\nu-C)}}\,.
\end{equation}
\new{In other words, the error bound is exponentially decaying if \mbox{$\nu>C \approx 3.69$}.}
\end{Theorem}

\emph{Proof}. Since the imaginary part of the integrand $\frac{1}{2}\,{\mathrm e}^{- \pi {\mathrm i} N t x}$, $t\in [-1, \,1]$, is odd, it holds
\begin{equation}
\label{eq:integralClenshawCurtis}
	\mathrm{sinc}(N \pi x) = \frac{1}{2} \int_{-1}^1 {\mathrm e}^{- \pi {\mathrm i} N t x}\, {\mathrm d}t
	= \frac{1}{2} \int_{-1}^1 \cos(\pi N t x)\, {\mathrm d}t\,.
\end{equation}
Therefore, we apply the Clenshaw--Curtis quadrature to the analytic function
{$f(t, x) \coloneqq \frac{1}{2}\, \cos(\pi N t x)$}, \mbox{$t \in [-1, 1]$},
with fixed parameter \mbox{$x \in [-1,\, 1]$}.
Note that it holds
$$
\sum_{k=0}^n w_k \,{\mathrm e}^{- \pi {\mathrm i} N z_k x} = \sum_{k=0}^n w_k \,\cos ( \pi N z_k x) + 0
$$
by the symmetry properties of the Chebyshev points $z_k$ and the coefficients $w_k$, namely \mbox{$z_k = - z_{n-k}$} and \mbox{$w_k =w_{n-k}$}, \mbox{$k=0,\ldots, n$} (see \cite[p.~359]{PlPoStTa18}).

By $E_{\rho}$ with some $\rho > 1$, we denote the Bernstein ellipse defined by
$$
E_{\rho} \coloneqq \bigg\{z \in \mathbb C:\,\mathrm{Re}\,z = \frac{1}{2}\bigg(\!\rho + \frac{1}{\rho}\!\bigg)\,\cos t\,, \;\mathrm{Im}\,z = \frac{1}{2}\bigg(\!\rho - \frac{1}{\rho}\!\bigg)\,\sin t\,, \;t\in [0,\,2\pi) \bigg\}\,.
$$
Then $E_{\rho}$ has the foci $-1$ and $1$. For simplicity, we choose $\rho = {\mathrm e}$.

For $z \in \mathbb C$ and fixed $x \in [-1,\,1]$, it holds
$$
\bigg| \frac{1}{2}\, \cos (\pi N x z)\bigg| \le \frac{1}{2}\, \cosh (\pi N x \, {\mathrm{Im}}\,z)\,.
$$
For  $z \in \mathbb C$ with ${\mathrm{Re}}\,z = 0$ we have
$$
\bigg| \frac{1}{2}\, \cos (\pi N x z)\bigg| = \frac{1}{2}\, \cosh (\pi N x \, {\mathrm{Im}}\,z)\,.
$$
Hence, in the interior of the Bernstein ellipse $E_{\mathrm e}$, the integrand is bounded, since
$$
\bigg| \frac{1}{2}\, \cos (\pi N x z)\bigg| \le \frac{1}{2}\, \cosh \frac{\pi N x \,({\mathrm e}^2 - 1)}{2\,\mathrm e} \le \frac{1}{2}\, \cosh \frac{\pi N \,({\mathrm e}^2 - 1)}{2\,\mathrm e}\,.
$$
\new{Therefore, by \cite[p.~146]{Tref13} we obtain the error estimate
\begin{equation}
\label{eq:estimateClenshawCurtis}
	\bigg|\mathrm{sinc}(N \pi x) - \sum_{k=0}^n w_k \,{\mathrm e}^{- \pi {\mathrm i} N z_k x}\bigg| \le \frac{144}{70\,({\mathrm e}^2 -1)}\, {\mathrm e}^{-n}\, \cosh \frac{\pi\,({\mathrm e}^2-1)\, N}{2\,{\mathrm e}}\,.
\end{equation}
By defining $C \coloneqq \frac{\pi\,({\mathrm e}^2-1)}{2\,{\mathrm e}}$,
the term $\mathrm e^{-n} \cosh(CN)$ in \eqref{eq:estimateClenshawCurtis} can be rewritten as
$$\mathrm e^{-n} \cosh(CN) = \mathrm e^{-\nu N} \cdot \tfrac 12(\mathrm e^{CN}+\mathrm e^{-CN})
= \tfrac 12 \,\mathrm e^{-N(\nu-C)} (1+\mathrm e^{-2CN})\,.$$
Thus, we end up with \eqref{eq:estimateClenshawCurtis2}.}
This completes the proof. \qedsymbol
\medskip

In practice, the coefficients $w_k$ in \eqref{eq:wk=} can be computed by a fast algorithm, the \mbox{so-called} discrete cosine transform of type~I (DCT--I) of length $n+1$, $n=2^t$, (see \cite[Algorithm~6.28 or Algorithm~6.35]{PlPoStTa18}).
This DCT--I uses the orthogonal cosine matrix of type~$\mathrm I$
$$
{\mathbf C}_{n+1}^{\mathrm I} \coloneqq \sqrt{\frac{2}{n}}\,\bigg(\varepsilon_n(j)\,\varepsilon_n(k)\,\cos \frac{jk \pi}{n}\bigg)_{j,k=0}^n\,.
$$

\begin{algorithm}[Fast computation of the coefficients $w_k$]\phantom{.}
	\label{alg:CompWeights}
	
	{\textit Input}:
	$n = 2^t$ with $t \in {\mathbb N} \setminus \{1\}$,
	$\varepsilon_n(0) = \varepsilon_n(n) \coloneqq \frac{\sqrt 2}{2}$, $\varepsilon_n(j) \coloneqq 1$ for \mbox{$j = 1, \ldots, n-1$}.
	\medskip
	
	$1$. Form the vector $(a_j)_{j=0}^n$ with $a_{2j} \coloneqq \varepsilon_n(2 j)\, \frac{2}{1 - 4 j^2}$, $j = 0, \ldots, n/2$ and \mbox{$a_{2j+1} \coloneqq 0$}, \mbox{$j= 0, \ldots, n/2 - 1$}. \\
	$2$.  Compute $({\hat a}_k)_{k=0}^n = {\mathbf C}_{n+1}^{\mathrm I} (a_j)_{j=0}^n$ by means of $\mathrm{DCT-I}$.\\
	$3$. Form the values $w_k \coloneqq \frac{1}{\sqrt{2 n}} \varepsilon_n(k) \, {\hat a}_k$, $k=0,\ldots, n$.
	\smallskip
	
	{\textit Output}: $w_k$ in \eqref{eq:wk=} for $k=0,\ldots,n$.
\end{algorithm}

\new{A similar approach can be found in }\cite{GrLeIn06}, where a Gauss--Legendre quadrature was applied to obtain explicit coefficients $w_k$ for given Legendre points $z_k$.
However, the computation of the coefficients $w_k$ using Algorithm~\ref{alg:CompWeights} is more effective for large $n$.

\begin{Example}
Now we visualize the result of Theorem~\ref{Thm:errorapproxsinc}.
\new{In Figure~\ref{fig:maxerr_approx_sinc_clenshawcurtis}~(a) the error bound \eqref{eq:estimateClenshawCurtis2} is depicted as a function of $N$ for several choices of \mbox{$\nu\in\{1,\dots,5\}$}, where $n = \nu N$.
It clearly demonstrates that \mbox{$\nu\geq 4$} is needed to obtain reasonable error bounds.}

Additionally, we compare the error constant and the maximum approximation error, cf.~\eqref{eq:estimateClenshawCurtis2}. 	
To measure the accuracy we consider a fine evaluation grid \mbox{$x_r=\frac{2r}{R}$}, \mbox{$r\in \I_R$}, with \mbox{$R\gg N$}, where \mbox{$R=3\cdot 10^5$} is fixed.
On this grid we calculate the discrete maximum error
\begin{equation}
\label{eq:err_approx_sinc_clenshawcurtis}
	\max_{r\in \I_R} \bigg| {\mathrm{sinc}}(\pi N x_r)- \sum_{k=0}^n w_k\,{\mathrm e}^{-\pi {\mathrm i}N z_k x_r} \bigg|\,
\end{equation}
for different bandwidths \mbox{$N=2^{\ell}$}, \mbox{$\ell=3,\dots,7$}.
For the parameter $n = \nu N$ we investigate several choices \mbox{$\nu\in\{1,\dots,10\}$}.	
We compute the coefficients $w_k$ using Algorithm~\ref{alg:CompWeights}.
Subsequently, the approximation to the $\mathrm{sinc}$ function is computed by means of the NFFT, which is possible since the $x_r$ are equispaced.
\new{
The results for both, the error bound \eqref{eq:estimateClenshawCurtis2} and the
maximum error \eqref{eq:err_approx_sinc_clenshawcurtis}, are displayed in Figure~\ref{fig:maxerr_approx_sinc_clenshawcurtis}~(b).}
It becomes apparent that for increasing oversampling factor~$\nu$, \new{the maximum error \eqref{eq:err_approx_sinc_clenshawcurtis} decreases to machine precision for all choices of $N$}.
\new{Even for rather large choices of $\nu$ (up to 10) the error remains stable, so there is no worsening in terms of $\nu$}.
\begin{figure}[t]
	\centering
	\captionsetup[subfigure]{justification=centering}
	\begin{subfigure}[t]{0.49\textwidth}
		\includegraphics[width=0.84\textwidth]{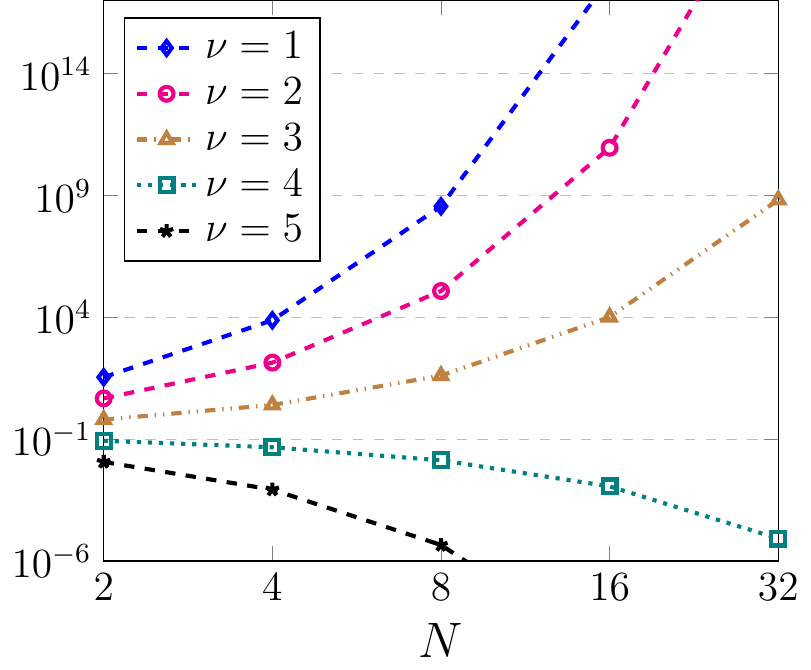}
		\caption{\new{with respect to the bandwidth $N$}}
	\end{subfigure}
	\begin{subfigure}[t]{0.49\textwidth}
		\includegraphics[width=0.87\textwidth]{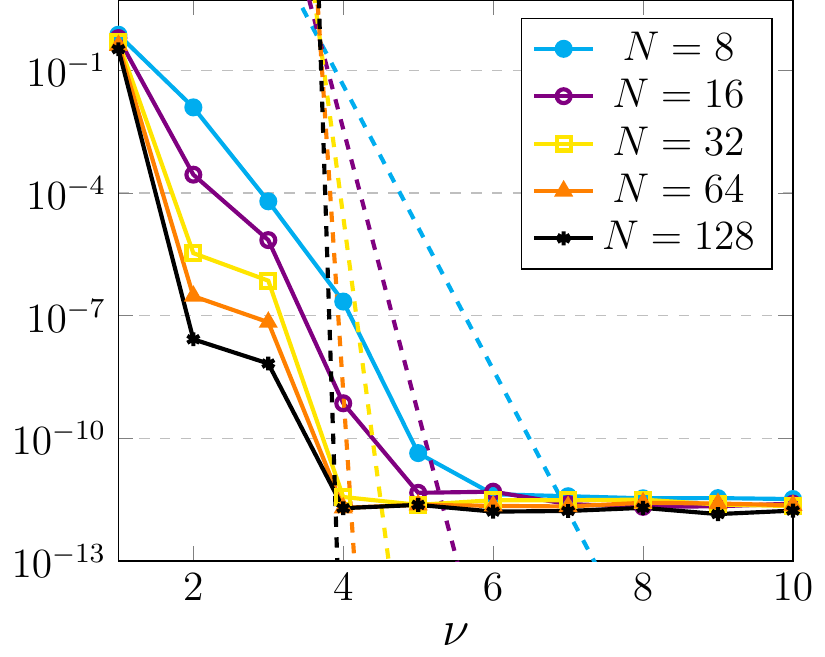}
		\caption{\new{with respect to the oversampling $\nu$}}
	\end{subfigure}
	\caption{\new{Error constant \eqref{eq:estimateClenshawCurtis2} (dashed) and maximum error \eqref{eq:err_approx_sinc_clenshawcurtis} (solid)} of the approximation of \mbox{$\mathrm{sinc}(N\pi x)$}, \mbox{$x\in[-1,1]$}, for different bandwidths \mbox{$N=2^{\ell}$}, \mbox{$\ell=1,\dots,7$}, where \mbox{$n=\nu N$}, \mbox{$\nu\in\{1,\dots,10\}$}, and Chebyshev nodes \mbox{$z_k\in[-1,1],\, k=0,\dots,n$}.
		\label{fig:maxerr_approx_sinc_clenshawcurtis}}
\end{figure}
\end{Example}
\newpage\vspace*{-2em}

\section{Discrete sinc transform} \label{Sec:discreteSinc}

Finally, we present an interesting signal processing application of the NNFFT.
If a signal \mbox{$h:\,\big[- \frac{1}{2},\, \frac{1}{2}\big] \to \mathbb C$} is to be reconstructed from its equispaced/nonequispaced
samples at \mbox{$a_k\in \big[- \frac{1}{2},\, \frac{1}{2}\big]$}, then $h$ is often modeled as linear combination of shifted $\mathrm{sinc}$ functions
\begin{equation}
\label{eq:h(x)=sum}
h(x) = \sum_{k \in \I_{L_1}} c_k \, \mathrm{sinc}\big(N \pi\,(x-a_k)\big)\,, \quad x \in \mathbb R\,,
\end{equation}
with complex coefficients $c_k$.
In the following, we propose a fast algorithm for the approximate computation of the \emph{discrete} $\mathrm{sinc}$ \emph{transform} (see \cite{GrLeIn06, LiBr11})
\begin{equation}
\label{eq:discretesinctransf}
h(b_{\ell}) = \sum_{k \in \I_{L_1}} c_k \, \mathrm{sinc}\big(N \pi\,(b_{\ell} - a_k)\big)\,, \quad \ell \in \I_{L_2}\,,
\end{equation}
where \mbox{$b_{\ell} \in \big[- \frac{1}{2},\, \frac{1}{2}\big]$} can be equispaced/nonequispaced points.

Such a function \eqref{eq:h(x)=sum} occurs by the application of the famous sampling theorem of Shannon--Whittaker--Kotelnikov (see e.\,g. \cite[pp.~86--88]{PlPoStTa18}).
Let \mbox{$f\in L_1(\mathbb R) \cap C(\mathbb R)$} be bandlimited on \new{\mbox{$\big[-\frac{L_2}{2},\,\frac{L_2}{2}\big]$}} for some \mbox{$L_2 > 0$}, i.\,e., the Fourier transform of $f$ is supported on \new{\mbox{$\big[-\frac{L_2}{2},\,\frac{L_2}{2}\big]$}}.
Then for \mbox{$N \in \mathbb N$} with \mbox{$N\ge L_2$}, the function $f$ is completely determined by its values \mbox{$f\big(\frac{k}{N}\big)$}, \mbox{$k \in \mathbb Z$}, and further $f$ can be represented in the form
$$
f(x) = \sum_{k\in \mathbb Z} f\big(\tfrac{k}{N}\big)\,\mathrm{sinc}\big(N \pi\big(x - \tfrac{k}{N}\big)\big)\,, \quad x \in \mathbb R\,,
$$
where the series converges absolutely and uniformly on $\mathbb R$.
By truncation of this series, we obtain the linear combination of shifted $\mathrm{sinc}$ functions
$$
\sum_{k\in \I_{L_1}} f\big(\tfrac{k}{N}\big)\,\mathrm{sinc}\big(N \pi\big(x - \tfrac{k}{N}\big)\big)\,, \quad x \in \mathbb R\,,
$$
which has the same form as \eqref{eq:h(x)=sum}, when $a_k$ are equispaced.

Since the naive computation of \eqref{eq:discretesinctransf} requires \mbox{$\mathcal O(L_1\cdot L_2)$} arithmetic operations, the aim is to find a more efficient method for the evaluation of \eqref{eq:discretesinctransf}.
Up to now, several approaches for a fast computation of the discrete $\mathrm{sinc}$ transform \eqref{eq:discretesinctransf} are known.
In \cite{GrLeIn06}, the discrete $\mathrm{sinc}$ transform \eqref{eq:discretesinctransf} is realized by applying a Gauss--Legendre quadrature rule to the integral~\eqref{eq:integralClenshawCurtis}. The result can then be approximated by means of two NNFFT's with \mbox{$\mathcal O((L_1+L_2)\log(L_1+L_2))$} arithmetic operations.
A multilevel algorithm with \mbox{$\mathcal{O}(L_2\log(1/\delta))$} arithmetic operations is presented in \cite{LiBr11} which is most effective for equispaced points $a_k$ and $b_\ell$ and, as the authors claim themselves, is only practical for rather large target evaluation accuracy \mbox{$\delta > 0$}.

In the following, we present a new approach for a fast $\mathrm{sinc}$ transform \eqref{eq:discretesinctransf}, where we approximate the function \mbox{$\mathrm{sinc}(N \pi x)$} by an exponential sum on the interval \mbox{$[-1, 1]$} by means of the Clenshaw--Curtis quadrature as described in Section~\ref{Sec:Approxsinc}.
Let the Chebyshev points \mbox{$z_j = \cos \frac{j \pi}{n}$}, \mbox{$j = 0,\ldots,n$},
and the coefficients $w_j$ defined by \eqref{eq:wk=} be given.
Utilizing \eqref{eq:estimateClenshawCurtis}, for arbitrary
\mbox{$a_k$, $b_{\ell} \in \big[- \frac{1}{2},\,\frac{1}{2}\big]$}
we obtain the approximation
\begin{equation*}
	\mathrm{sinc}\big(N \pi(a_k - b_\ell)\big)
	\approx
	\sum_{j=0}^n w_j\, \,\mathrm e^{-\pi\mathrm i N z_j (a_k - b_\ell)}
	=
	\sum_{j=0}^n w_j \,\mathrm e^{-\pi\mathrm i N z_j a_k}\,\mathrm e^{\pi\mathrm i N z_j b_\ell}\,.
\end{equation*}
Inserting this approximation into \eqref{eq:discretesinctransf} yields
\begin{align}
\label{eq:approxsinctrafo}
h_{\ell} &\coloneqq \sum_{k\in \I_{L_1}} c_k\, \sum_{j=0}^n w_j \,\mathrm e^{-\pi\mathrm i N z_j a_k}\,\mathrm e^{\pi\mathrm i N z_j b_\ell} \notag\\
	&=
	\sum_{j=0}^n w_j
	\bigg(\sum_{k\in \I_{L_1}} c_k \,\mathrm e^{-\pi\mathrm i N z_j a_k}\bigg)
	\,\mathrm e^{\pi\mathrm i N z_j b_\ell}\,,
	\quad \ell\in \I_{L_2}\,.
\end{align}
If $\varepsilon > 0$ denotes a target accuracy, then we choose $n = 2^t$, $t \in \mathbb N \setminus \{1\}$ such that by Theorem~\ref{Thm:errorapproxsinc} it holds
$$
\new{\frac{36\,(1+\mathrm e^{-2CN})}{35\,({\mathrm e}^2 -1)} \,\mathrm e^{-N(\nu-C)}} < \varepsilon\,, \quad \nu > C = 3.69\,.
$$
For example, in the case $\varepsilon = 10^{-8}$ we obtain $n \ge 4N$ for $N\geq 54$.

We recognize that the term inside the brackets of \eqref{eq:approxsinctrafo} is an exponential sum of the form \eqref{eq:f(xj)}, which can be computed by means of an NNFFT.
Then the resulting outer sum is of the same form such that this can also be computed by means of an NNFFT.
Thus, as in \cite{GrLeIn06} we may compute the discrete $\mathrm{sinc}$ transform \eqref{eq:discretesinctransf} by means of an NNFFT, a multiplication by the precomputed coefficients $w_j$ as well as another NNFFT afterwards.
Hence, the fast $\mathrm{sinc}$ transform, which is an application of the NNFFT, can be summarized as follows.

\medskip
\begin{algorithm}[Fast $\boldsymbol{\mathrm{sinc}}$ transform]\phantom{.}
\label{alg:fastsinc}

	{\textit Input}: \mbox{$N \in \mathbb N$}, \mbox{$L_1$, $L_2\in 2{\mathbb N}$}
	as well as
	\mbox{$c_k\in \mathbb{C}$}, \mbox{$a_k \in \big[-\frac{1}{2},\,\frac{1}{2}\big]$}
	for \mbox{$k\in \I_{L_1}$},
	\mbox{$z_j =\cos \frac{j \pi}{n}$} with \mbox{$j =0,\ldots,n$} and \mbox{$n\geq 4N$}.
	\medskip
	
	$0$. Precompute the values $w_j$, $j =0,\ldots,n$, by Algorithm~$\mathrm{\ref{alg:CompWeights}}$.
	\smallskip
	
	$1$. For all $j=0,\dots,n,$ compute by $\mathrm{NNFFT}$
	\begin{equation*}
		 g_j \coloneqq \sum_{k\in \I_{L_1}} c_k \,\mathrm e^{-\pi\mathrm i N z_j a_k}\,,
	\end{equation*}
    where ${\tilde g}_j$ is the approximate value of $g_j$.\\
	$2$. For all $j=0,\dots,n,$ form the products
	\begin{equation*}
		\alpha_j \coloneqq w_j \cdot {\tilde g}_{j}\,.
	\end{equation*}
	$3$. For all $\ell\in \I_{L_2}$ compute by $\mathrm{NNFFT}$
	\begin{equation}
	\label{eq:output_sinc_trafo}
		{\hat h}_\ell \coloneqq \sum_{j=0}^n \alpha_j \,\mathrm e^{ \pi\mathrm i N z_j b_\ell}\,,
	\end{equation}
    where ${\tilde h}_{\ell}$ is the approximate value of ${\hat h}_\ell$.
	\smallskip
	
	{\textit Output}: $\tilde h_\ell$ approximate value of \eqref{eq:discretesinctransf} for $\ell\in \I_{L_2}$.
\end{algorithm}

If we use the same NNFFT in both steps (with the window functions $\varphi_j$, truncation parameters $m_j$, and oversampling factors $\sigma_j$ for \mbox{$j =1,\,2$}),  Algorithm~\ref{alg:fastsinc} requires
all in all
$$
\mathcal O(N\log N + L_1+L_2+2n)
$$
arithmetic operations.

Considering the discrete $\mathrm{sinc}$ transform \eqref{eq:discretesinctransf}, we can deal with the special sums of the form
\begin{equation*}
	h\big(\tfrac{\ell}{N}\big)= \sum_{k\in \I_{L_1}} c_k \, \mathrm{sinc} \big(N\pi \big(a_k-\tfrac{\ell}{N}\big)\big)\,,\quad
	\ell\in \I_{N},
\end{equation*}
i.\,e., we are given equispaced points \mbox{$b_\ell = \frac{\ell}{N}$} with \mbox{$L_2=N$}. In this special case, we simply obtain an adjoint NFFT instead of the NNFFT in step~3 of Algorithm~\ref{alg:fastsinc}.
Therefore, the computational cost of Algorithm~\ref{alg:fastsinc} reduces to \mbox{$\mathcal O(N\log N + L_1 + n)$}.
In the case, where \mbox{$a_k = \frac{k}{L_1}$}, \mbox{$k\in \I_{L_1}$}, the NNFFT in step~1 of Algorithm~\ref{alg:fastsinc} naturally turns into an NFFT. Clearly, in this case the same amount of arithmetic operations is needed as in the first special case.
If both sets of nodes $a_k$ and $b_\ell$ are equispaced, then the computational cost reduces even more to \mbox{$\mathcal O(N\log N + n)$}.
Hence, these modifications are automatically be included in our fast $\mathrm{sinc}$ transform.

A quite similar approach was already developed in \cite{AlAu15} for the computation of the Coulombian interaction between punctual masses, where the main idea is using two different quadrature rules to approximate the given problem.
Then the computation can be done by means of NNFFTs, i.\,e., they receive a 3-step method analogous to Algorithm~\ref{alg:fastsinc}.

Now we study the error of the fast $\mathrm{sinc}$ transform in Algorithm~\ref{alg:fastsinc}, which is measured in the form	
\begin{equation}
\label{eq:errsinctrafo}
	\max_{\ell \in \I_{L_2}}\ | h(b_\ell)-\tilde h_\ell |\,.
\end{equation}

\begin{Theorem}
Let \mbox{$N\in \mathbb N$} with \mbox{$N \gg 1$} and \mbox{$L_1$, $L_2 \in 2 \mathbb N$} be given.
Let \mbox{$N_1 = \sigma_1 N \in 2 \mathbb N$} with \mbox{$\sigma_1 > 1$}.
For fixed \mbox{$m_1 \in \mathbb N \setminus \{1\}$} with \mbox{$2 m_1 \ll N_1$},
let  \mbox{$N_2 = \sigma_2\,(N_1 + 2m_1)$} with \mbox{$\sigma_2 > 1$}.
For \mbox{$m_2 \in \mathbb N \setminus \{1\}$} with
\mbox{$2m_2 \le \big(1 - \frac{1}{\sigma_1}\big)\, N_2$}, let
$\varphi_1$ and $\varphi_2$ be the window functions of the form \eqref{eq:twowindows}.
Let \mbox{$a_k$, $b_{\ell} \in \big[ - \frac{1}{2},\,\frac{1}{2}\big]$}
with \mbox{$k \in \I_{L_1}$}, \mbox{$\ell \in \I_{L_2}$} be arbitrary points
and let \mbox{$c_k \in \mathbb C$}, \mbox{$k \in \I_{L_1}$}, be arbitrary coefficients.
Let \mbox{$a > 1$} be the constant \eqref{eq:consta}.
For a given target accuracy \mbox{$\varepsilon > 0$}, the number
\mbox{$n = 2^t$}, \mbox{$t \in \mathbb N \setminus \{1\}$}, is chosen such that
\begin{equation}
\label{eq:condvarepsilon}
\new{\frac{36\,(1+\mathrm e^{-2CN})}{35\,({\mathrm e}^2 -1)} \,\mathrm e^{-N(\nu-C)}} < \varepsilon\,, \quad \new{\nu > C = 3.69}\,.
\end{equation}

Then the error of the fast $\mathrm{sinc}$ transform can be estimated by
\begin{eqnarray}
\label{eq:errorsinctransf}
\max_{\ell \in \I_{L_2}} \big| h(b_\ell)- {\tilde h}_\ell \big| &\le& \Bigg(\varepsilon + 2 \bigg[E_{\sigma_1}(\varphi_1) + \frac{a}{ {\hat \varphi}_1\big(\frac{N}{2}\big)}\,E_{\sigma_2}(\varphi_2)\bigg] \nonumber\\
& & + \bigg[E_{\sigma_1}(\varphi_1) + \frac{a}{ {\hat \varphi}_1\big(\frac{N}{2}\big)}\,E_{\sigma_2}(\varphi_2)\bigg]^2 \Bigg)\, \sum_{k\in \I_{L_1}} |c_k|\,,
\end{eqnarray}
where $E_{\sigma_j}(\varphi_j)$ for $j = 1,\,2$, are the general $C(\mathbb T)$-error constants of the form \eqref{eq:Esigma(varphi)}. If it holds
\vspace*{-2ex}
\begin{align}
\label{eq:assump}
	E_{\sigma_1}(\varphi_1) + \frac{a}{ {\hat \varphi}_1\big(\frac{N}{2}\big)}\,E_{\sigma_2}(\varphi_2) \le 1\,,
\end{align}
one can use the simplified estimate
\begin{align}
\label{eq:errorsinctransf_simplified}
	\max_{\ell \in \I_{L_2}} \big| h(b_\ell)- {\tilde h}_\ell \big| \le \bigg(\varepsilon + \new{3} E_{\sigma_1}(\varphi_1) + \frac{\new{3}a}{ {\hat \varphi}_1\big(\frac{N}{2}\big)}\,E_{\sigma_2}(\varphi_2)\bigg) \,  \sum_{k\in \I_{L_1}} |c_k|\,.
\end{align}
\end{Theorem}

\emph{Proof.} By \eqref{eq:approxsinctrafo}, the value $h_{\ell}$ is an approximation of $h(b_{\ell})$.
Since \mbox{$a_k$, $b_{\ell}\in \big[- \frac{1}{2}, \,\frac{1}{2}\big]$},
it holds by \eqref{eq:estimateClenshawCurtis} and \eqref{eq:condvarepsilon} that
$$
\bigg| \mathrm{sinc}\big(\pi N (a_k - b_{\ell})\big) - \sum_{j=0}^n w_j\, {\mathrm e}^{-\pi {\mathrm i}Nz_j (a_k - b_{\ell})}\bigg| \le \varepsilon\,.
$$
Hence, we conclude that
\begin{equation}
\label{eq:estim1}
|h(b_{\ell}) - h_{\ell}| \le \varepsilon \sum_{k \in \I_{L_1}} |c_k|\,, \quad \ell \in \I_{L_2}\,.
\end{equation}
After step~1 of Algorithm~\ref{alg:fastsinc}, the error of the NNFFT (with the window functions $\varphi_1$ and $\varphi_2$) can be estimated by Theorem~\ref{Thm:errorestimate} in the form
$$
|g_j - {\tilde g}_j| \le \bigg[E_{\sigma_1}(\varphi_1) + \frac{a}{ {\hat \varphi}_1\big(\frac{N}{2}\big)}\,E_{\sigma_2}(\varphi_2)\bigg] \sum_{k \in \I_{L_1}} |c_k|\,, \quad j = 0,\ldots, n\,.
$$
Using \eqref{eq:sumwk=1}, step~2 of Algorithm~\ref{alg:fastsinc} generates the error
\begin{eqnarray}
\label{eq:estim2}
|{\hat h}_{\ell} - h_{\ell} | &\le& \sum_{j=0}^n w_j \,|g_j - {\tilde g}_j| \nonumber 
\le \bigg(\sum_{j=0}^n w_j\bigg) \bigg[E_{\sigma_1}(\varphi_1) + \frac{a}{ {\hat \varphi}_1\big(\frac{N}{2}\big)}\,E_{\sigma_2}(\varphi_2)\bigg] \sum_{k \in \I_{L_1}} |c_k| \nonumber \\
&=& \bigg[E_{\sigma_1}(\varphi_1) + \frac{a}{ {\hat \varphi}_1\big(\frac{N}{2}\big)}\,E_{\sigma_2}(\varphi_2)\bigg] \sum_{k \in \I_{L_1}} |c_k|\,.
\end{eqnarray}
After step~3 of Algorithm~\ref{alg:fastsinc}, the error of the NNFFT (with the same window functions $\varphi_1$ and $\varphi_2$) can be estimated by Theorem~\ref{Thm:errorestimate} in the form
$$
\big|{\hat h}_{\ell} - {\tilde h}_{\ell}\big| \le \bigg[E_{\sigma_1}(\varphi_1) + \frac{a}{ {\hat \varphi}_1\big(\frac{N}{2}\big)}\,E_{\sigma_2}(\varphi_2)\bigg] \sum_{j=0}^n w_j \big|{\tilde g}_j\big|\,, \quad \ell \in \I_{L_2}\,.
$$
Using the triangle inequality, we obtain
\begin{eqnarray*}
\big|{\tilde g}_j\big| &\le& |g_j| + |g_j - {\tilde g}_j| \le \sum_{k \in \I_{L_1}} |c_k| + |g_j - {\tilde g}_j|\\
&\le& \sum_{k \in \I_{L_1}} |c_k| + \bigg[E_{\sigma_1}(\varphi_1) + \frac{a}{ {\hat \varphi}_1\big(\frac{N}{2}\big)}\,E_{\sigma_2}(\varphi_2)\bigg]\sum_{k \in \I_{L_1}} |c_k|\,, \quad j=0,\ldots,n
\end{eqnarray*}
such that by \eqref{eq:sumwk=1}
\begin{eqnarray}
\label{eq:estim3}
\big|{\hat h}_{\ell} - {\tilde h}_{\ell}\big| &\le& \bigg[E_{\sigma_1}(\varphi_1) + \frac{a}{ {\hat \varphi}_1\big(\frac{N}{2}\big)}\,E_{\sigma_2}(\varphi_2)\bigg]\bigg(\sum_{j=0}^n w_j\bigg)\sum_{k \in \I_{L_1}} |c_k| \nonumber\\
& & +\, \bigg[E_{\sigma_1}(\varphi_1) + \frac{a}{ {\hat \varphi}_1\big(\frac{N}{2}\big)}\,E_{\sigma_2}(\varphi_2)\bigg]^2 \bigg(\sum_{j=0}^n w_j\bigg)\sum_{k \in \I_{L_1}} |c_k| \nonumber\\
&=& \bigg[E_{\sigma_1}(\varphi_1) + \frac{a}{ {\hat \varphi}_1\big(\frac{N}{2}\big)}\,E_{\sigma_2}(\varphi_2)\bigg]\,\sum_{k \in \I_{L_1}} |c_k| \nonumber\\
& & +\, \bigg[E_{\sigma_1}(\varphi_1) + \frac{a}{ {\hat \varphi}_1\big(\frac{N}{2}\big)}\,E_{\sigma_2}(\varphi_2)\bigg]^2 \,\sum_{k \in \I_{L_1}} |c_k| \,.
\end{eqnarray}
Thus, the error of Algorithm~\ref{alg:fastsinc} can be estimated by
$$
\big| h(b_{\ell}) - {\tilde h}_{\ell}\big| \le |h(b_{\ell}) - h_{\ell}| + |h_{\ell} - {\hat h}_{\ell}| + \big|{\hat h}_{\ell} - {\tilde h}_{\ell}\big|\,, \quad \ell \in \I_{L_2}\,.
$$
From \eqref{eq:estim1} -- \eqref{eq:estim3} it follows the estimate \eqref{eq:errorsinctransf}. 
\new{If it holds \eqref{eq:assump}, we have
\begin{eqnarray*}
	\bigg[E_{\sigma_1}(\varphi_1) + \frac{a}{ {\hat \varphi}_1\big(\frac{N}{2}\big)}\,E_{\sigma_2}(\varphi_2)\bigg]^2 
	&\le& E_{\sigma_1}(\varphi_1) + \frac{a}{ {\hat \varphi}_1\big(\frac{N}{2}\big)}\,E_{\sigma_2}(\varphi_2)
\end{eqnarray*}
and therefore the simplified estimate \eqref{eq:errorsinctransf_simplified}.}
\qedsymbol
\medskip

\new{Thus, the error of Algorithm~\ref{alg:fastsinc} for the fast $\mathrm{sinc}$ transform mostly depends on the target accuracy $\varepsilon$ of the precomputation and on the general $C(\mathbb T)$-error constants $E_{\sigma_j}(\varphi_j)$,
$j = 1,\,2$, of the window functions $\varphi_j$, $j = 1,\,2$, see Theorem \ref{Thm:errorestimate}.}

\begin{Example}
Next we verify the accuracy of our fast $\mathrm{sinc}$ transform in Algorithm~\ref{alg:fastsinc}.
To this end, we choose random nodes \mbox{$a_k\in \big[-\frac 12,\frac 12\big]$}, equispaced points \mbox{$b_\ell = \frac{\ell}{N}$} with \mbox{$\ell \in \I_N$}, as well as random coefficients \mbox{$c_k\in\mathbb C$}, \mbox{$k\in \I_{L_1}$}, and compute
the discrete $\mathrm{sinc}$ transform \eqref{eq:discretesinctransf} directly as well as its approximation \eqref{eq:output_sinc_trafo} by means of the fast $\mathrm{sinc}$ transform.
Subsequently, we compute the maximum error \eqref{eq:errsinctrafo}.
Due to the randomness of the given values this test is repeated one hundred times and afterwards the maximum error over all repetitions is computed.

In this experiment we choose different bandwidths \mbox{$N=2^k$}, \mbox{$k=5,\dots,13,$}
and without loss of generality we use \mbox{$L_1=\frac N2$}.
We apply Algorithm~\ref{alg:fastsinc} using the weights $w_j$ computed by means of Algorithm~\ref{alg:CompWeights} and the Chebyshev points \mbox{$z_j =  \cos \frac{j \pi}{n}$}, \mbox{$j = 0,\ldots, n$}.
Therefore, we only have to examine the parameter choice of \mbox{$n\ge 4N$}.
To this end, we compare the results for several choices, namely for \mbox{$n\in\{4N,6N,8N\}$}.
The appropriate results can be found in Figure~\ref{fig:maxerr_sinc_trafo_clenshawcurtis}.
We see that for large $N$ there is almost no difference between the different choices of $n$.
However, we point out that a higher choice heavily increases the computational cost of Algorithm~\ref{alg:fastsinc}.
Therefore, it is recommended to use the smallest possible choice \mbox{$n=4N$}.
Compared to \cite{GrLeIn06} the same approximation errors are obtained, but with a more efficient precomputation of weights.
\begin{figure}[ht]
	\centering
	\captionsetup[subfigure]{justification=centering}
	\begin{subfigure}[t]{0.42\textwidth}
		\includegraphics[width=\textwidth]{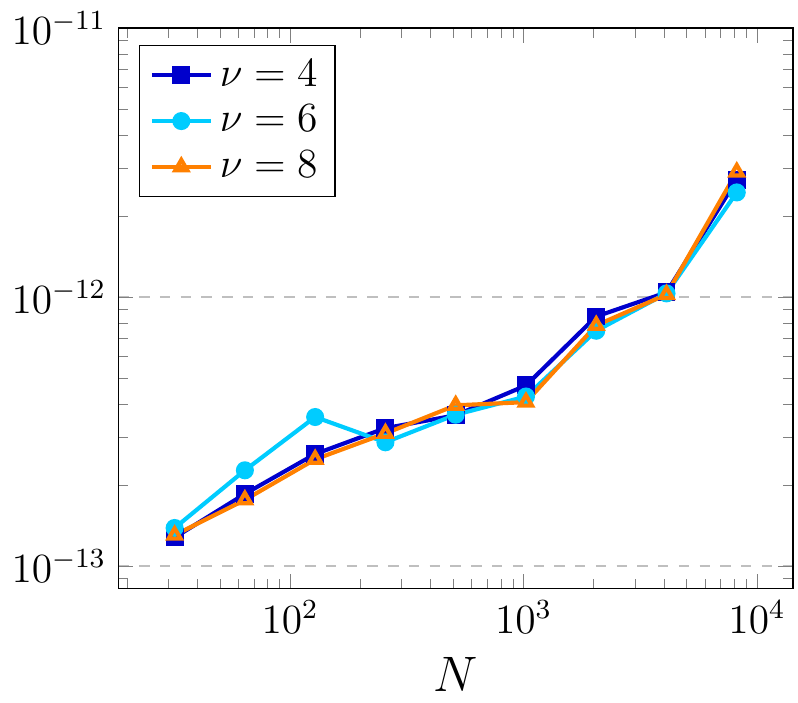}
	\end{subfigure}
	\caption{Maximum error \eqref{eq:errsinctrafo} for several bandwidths \mbox{$N=2^k$}, \mbox{$k=5,\dots,13,$}, shown for \mbox{$n=\nu N$}, \mbox{$\nu\in\{4,6,8\}$}, using the coefficients $w_j$ obtained by  Algorithm~\ref{alg:CompWeights}.
		\label{fig:maxerr_sinc_trafo_clenshawcurtis}}
\end{figure}
\end{Example}

\section*{Acknowledgments}
Melanie Kircheis gratefully acknowledges the funding support from the European Union and the Free State of Saxony (ESF).
Daniel Potts acknowledges funding by Deutsche Forschungsgemeinschaft (German Research Foundation) -- Project--ID 416228727 -- SFB 1410.

Moreover, the authors thank the referees and the editor for their very helpful suggestions for improvements.

\end{document}